\documentclass[10pt]{article}

\usepackage[latin1]{inputenc}
\usepackage[spanish,activeacute]{babel}
\usepackage[spanish]{layout}
\usepackage[reqno]{amsmath}
\usepackage{amsthm}
\usepackage{amssymb}
\usepackage{array}
\usepackage{longtable}
\usepackage{enumerate}
\usepackage{babelbib}
\usepackage[format=plain,labelsep=period,font={footnotesize},width=0.8\textwidth]{caption}
\usepackage{verbatim}
\usepackage{fancyhdr}
\usepackage{titlesec}
\usepackage{textcomp}
\usepackage{authblk}
\usepackage[papersize={8.5in,11in},left=1cm,right=1cm,top=1cm,bottom=1cm,includefoot=true,
ignorehead=false,ignoremp=true,driver=dvips,footskip=1cm,headheight=0cm,headsep=0cm]{geometry}
\usepackage[pdftex]{graphicx}
\usepackage[pdftex]{hyperref}
\hypersetup{bookmarksnumbered=true,colorlinks=true,linkcolor=red,citecolor=black,bookmarksopen=true,{pdftitle=Resoluci\'on de tri\'angulos oblicu\'angulos},pdfauthor=Diego Fernando Ram\'irez}

\titleformat{\section}[hang]{\normalsize\center\bfseries}{\thesection.}{6pt}{}
\titleformat{\subsection}[hang]{\normalsize\center\it}{\thesubsection}{6pt}{}
\titleformat{\subsubsection}[runin]{\normalsize\it}{\thesubsubsection}{2pt}{}

\providecommand{\al}{\alpha}
\newcommand{\GC}[1]{\biggl[#1\biggr]}   %Corchetes grandes
\newcommand{\NC}[1]{\bigl[#1\bigr]}     %Corchetes normales

\newcommand{\GP}[1]{\biggl(#1\biggr)}   %Parentesis grandes
\newcommand{\seg}[1]{\overline{#1}}

\allowdisplaybreaks
\date{}
\title{Resoluci\'on de tri\'angulos oblicu\'angulos}
\author{Diego Fernando Ram\'irez Jim\'enez\\{\tt df.ramirezj@uniandes.edu.co}}
\begin{document}
\maketitle
\thispagestyle{empty}
\pagestyle{fancy}
\fancyhead{}
\fancyfoot{}
\fancyfoot[CF]{\thepage}
\renewcommand{\headrulewidth}{0.0pt}
\renewcommand{\footrulewidth}{0.4pt}
%contadores
%\renewcommand{\thefigure}{\arabic{figure}}
\nocite{*}
%correccion numeros romanos
\makeatletter
\def\@roman#1{\romannumeral#1}
\makeatother
%%%%%%%%% Documento %%%%%%%%%%
\addcontentsline{toc}{section}{Resumen}
\textbf{Resumen:} Se enuncia los principales teoremas empleados en la resoluci'on de tri'angulos oblicu'angulos. Con ellos, se ilustra c'omo resolver los cinco casos de resoluci'on que se presentan, incluyendo algunos caso at'ipicos (cuando se conoce el per'imetro y dos 'angulos internos, o un lado, el 'angulo comprendido entre los lados restantes, as'i como su suma). Luego, se discute la determinaci'on del 'area para cada caso, las relaciones entre los radios de las circunferencias inscritas, circunscritas y excritas, las longitudes de las medianas, bisectrices y alturas. 

\emph{T'erminos claves:} Tri'angulo oblicu'angulo, 'angulo interno, lado, semiper'imetro, circunferencia inscrita, circunferencia circunscrita, circunferencia excrita, mediana, bisectriz, altura, tri'angulo excrito, tri'angulo pedal. 

\textbf{Abstract:} {The principal theorems for solving oblicue triangles are presented. We shall show how to solve the five classical cases, and also some atipical cases, for instance, known the perimeter and two internal angles; or a lade, an angle between the other lades as like their sum. Finally, we shall discute how to calcule the area of triangles for the difference cases, the relations between the radii of inscribed, circunscribed and excribes  circles, the lengths of medians, bisectors and heights.}

\emph{Key terms:} {Oblicue triangle, internal angle, lade, semiperimeter, inscribed circle, circumscribed circle, excrite circle, median, bisector, height.} 
\section{Introducci'on}
La rama de la trigonometr'ia conocida como la resoluci'on de tri'angulos es el espacio ideal para aplicar muchos de sus resultado, adem'as de permitir integrar otros campos de la matem'atica elemental como el 'algebra y la geometr'ia. Esto se traduce en una de las primeras oportunidades que tiene un estudiante que comienza a aprender trigonometr'ia de ver todo su potencial. En muchos textos sobre la materia, el tema en cuesti'on era tratado con una extensi'on considerable, en parte, porque los c'alculos que se emplean son bastante complicados y requer'ian el empleo de los logaritmos para llevarlos a cabo. Sin embargo, la llegada de nuevas tecnolog'ias hizo a un lado estos inconvenientes, simplificando la ense'nanza de este tema, pero a tal grado que muchos resultados como el teorema de la tangente, las f'ormulas de Mollweide, entre otras, dejaron de ense'narse, aún cuando hacen parte integral de este campo; y en los textos en donde se mencionan, est'an bastante entrelazados con el c'alculo logar'itmico, que si bien, fue una herramienta crucial en su momento, hace que su discusi'on resulte bastante engorrosa. Guiado por estas reflexiones, me di a la tarea de escribir estas l'ineas con el prop'osito de proporcionar a quienes quieran profundizar en el tema un material de referencia.

Con este objetivo, el presente art'iculo est'a dividido en seis secciones. En la primera, se define lo que se entiende por tri'angulos oblicu'angulos, se indica la notaci'on a emplear y los casos cl'asicos de resoluci'on. En la segunda, se enuncian y demuestran los teoremas fundamentales usados para resolver tri'angulos en general. Entre ellos est'an el del seno, del coseno, de la tangente, las f'ormulas de Mollweide y las de proyecci'on. En la siguiente secci'on, se resuelve los cinco casos de resoluci'on de tri'angulos. Las discusiones anal'iticas para cada caso se acompa'nan con sus correspondientes construcciones geom'etricas. En esta secci'on se trata, adem'as, de dos casos cuyos par'ametros dados no corresponden a los casos cl'asicos. En particular, se muestra c'omo resolver un tri'angulo cuando se conoce su per'imetro y dos de sus 'angulos internos, y cuando se conoce un lado, su 'angulo opuesto y la suma de los lados restante. La cuarta secci'on se enfoca en el c'alculo del 'area de un tri'angulo para cada caso tratado. La siguiente secci'on extiende la discusi'on de los radios de las diferentes circunferencias asociadas a un tri'angulo, tema que surge en las demostraciones de los teoremas enunciados. La 'ultima secci'on tiene como fin mostrar c'omo se calculan las medianas, bisectrices y alturas en un tri'angulo, adem'as de mostrar c'omo se resuelve un tri'angulo cuando se conocen las tres medianas o las tres alturas.
\section{Definiciones y notaci'on}\label{S1}
\begin{figure}[htb]
\centering\fbox{
\includegraphics[scale=0.8]{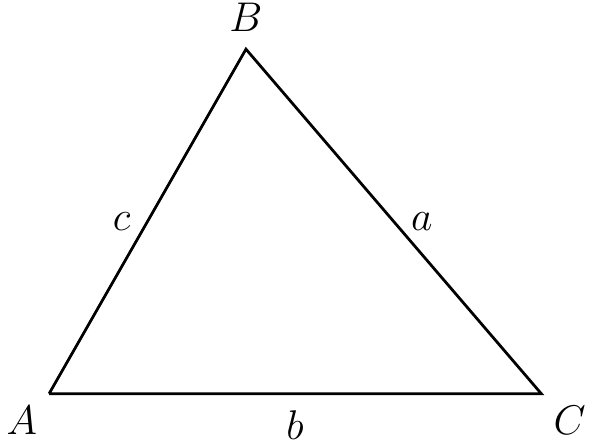}}
\caption{Ilustraci'on de la convenci'on para nombrar los lados y los v'ertices de un tri'angulo oblicu'angulo.}\label{F1}
\end{figure}

Se entiende por \emph{tri'angulo oblicu'angulo} por uno que sea acut'angulo u obtus'angulo\footnote{En este art'iculo no se discute la resoluci'on de tri'angulos rect'angulos, pues los problemas referentes a ellos se resuelven exitosamente usando razones trigonom'etricas. Adem'as, un tri'angulo oblicu'angulo puede resolverse con tri'angulos rect'angulos, sin la ayuda de los teoremas enunciados en la secci'on \ref{S2}. El lector interesado en estas metodolog'ias, puede ver, por ejemplo, \cite{B14}, p'ag. 110--114.}. Para denotar los lados y los v'ertices de un tri'angulo oblicu'angulo se sigue la siguiente convenci'on: para los lados, se emplean letras min'usculas y para los v'ertices de los 'angulos opuestos a estos lados se utilizan las correspondientes letras en may'uscula. En la figura \ref{F1}, el lector ver'a materializado el convenio en cuesti'on. En algunos tratados sobre este tema, el valor de los 'angulos internos suele indicarse mediante letras del alfabeto griego, pero aqu'i, dichos valores se representar'an con las letras asignadas para nombrar los v'ertices del tri'angulo. Si bien, esto puede parecer al lector un abuso de notaci'on, permite plantear las ecuaciones enunciadas de los teoremas referentes a la resoluci'on de tri'angulos oblicu'angulos de tal manera que, de una ecuaci'on referente a un teorema en particular, se pueda escribir las dem'as mediante un intercambio c'iclico de letras. 

Una vez establecida la nomenclatura, es menester discutir cu'ales son los par'ametros m'inimos que deben conocerse del tri'angulo para determinarlo completamente, y adem'as, si dicha elecci'on conduce a una, varias o a ninguna soluci'on. Por supuesto, el inter'es se centra en aquellos que proporcionen al menos una sola. Sobra advertir que una vez establecidos, no siempre se garantiza su unicidad. Bajo estas circunstancias, la geometr'ia plana dice que el n'umero de estos deben ser tres, de los cuales uno de ellos debe ser un lado. Dependiendo de los par'ametros conocidos se tendr'an cinco casos a saber: 
\begin{itemize}
	\item[-] \emph{Caso 1:} dos 'angulos y un lado com'un a ellos;
	\item[-] \emph{Caso 2:} dos 'angulos y un lado opuesto a uno de ellos;
	\item[-] \emph{Caso 3:} dos lados y un 'angulo opuesto a uno de ellos;
	\item[-] \emph{Caso 4:} dos lados y el 'angulo comprendido entre ellos; 
	\item[-] \emph{Caso 5:} tres lados.
\end{itemize}
Establecidos los casos a tratar, ahora la discusi'on se dirige a encontrar teoremas que permitan resolverlos. Los tres primeros casos se resuelven completamente mediante el teorema del seno, el cuarto se soluciona con el teorema del coseno o la combinaci'on de los teoremas de la tangente y del seno y el quinto se trata ya sea con el del coseno o con las f'ormulas del 'angulo medio.
\section{Teoremas referentes a tri'angulos oblicu'angulos}\label{S2}
A continuaci'on, se enuncian y demuestran los teoremas m'as importantes sobre tri'angulos oblicu'angulos.
\subsection{El teorema del seno}\label{S2-1}
El teorema del seno dice lo siguiente: \textit{La longitud de un lado cualquiera de un tri'angulo oblicu'angulo es proporcional al seno de su 'angulo opuesto}.

De acuerdo a la figura \ref{F1}, este teorema se escribe de la siguiente manera:
\begin{equation}\label{E1}
	\frac{a}{\sen{A}}=\frac{b}{\sen{B}}=\frac{c}{\sen{C}}
\end{equation}
\begin{figure}[htb]
\centering
\begin{tabular}{cc}
\fbox{\includegraphics[scale=0.8]{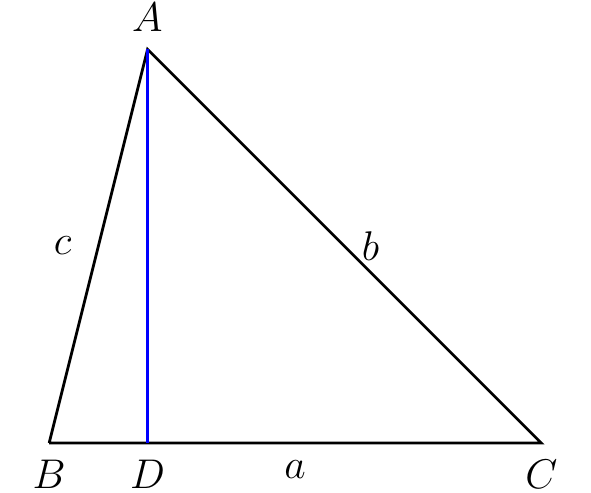}}&\fbox{\includegraphics[scale=0.8]{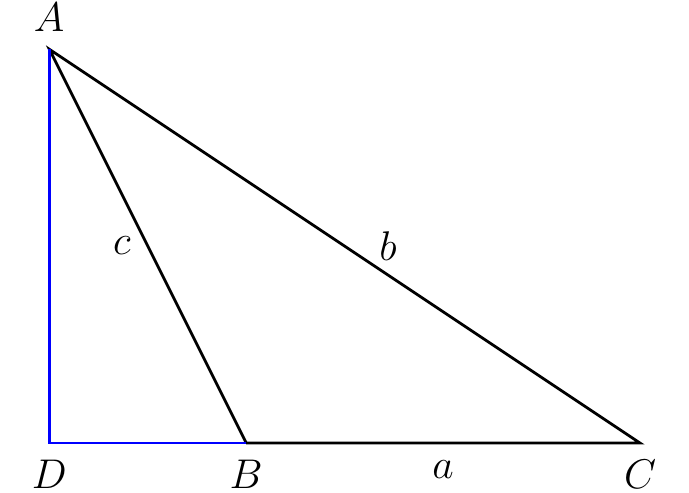}}\\
a) & b)
\end{tabular}
\caption{Construcci\'on geom\'etrica para la demostraci\'on de los teoremas del seno y del coseno.}\label{F2}
\end{figure}
Demostraci'on: La prueba de este teorema se divide en dos casos:

Caso I: El tri'angulo es acut'angulo. Consid'erese el tri'angulo acut'angulo de la figura \ref{F2}-a, al cual se le traza la altura correspondiente al lado $a$, la cual lo intercepta en $D$. Del tri'angulo rect'angulo $ABD$:
\[
\seg{AD}=c\,\sen{B},
\]
y del tri'angulo rect'angulo $ACD$:
\[
\seg{AD}=b\,\sen{C}.
\]
Igualando las dos expresiones anteriores, se tiene:
\[
b\,\sen{C}=c\,\sen{B}.
\]
Dividiendo por $\sen{B}\sen{C}$, se tiene, finalmente:
\[
\frac{b}{\sen{B}}=\frac{c}{\sen{C}}
\]
Caso II: El tri'angulo es obtus'angulo. En el tri'angulo obtus'angulo de la figura \ref{F2}-b, al cual se le traza la altura correspondiente al lado $a$, la cual lo intercepta en $D$.
Del tri'angulo rect'angulo $ABD$:
\[
\seg{AD}=c\,\sen{(180-B)}=c\,\sen{B},
\]
y del tri'angulo rect'angulo $ACD$:
\[
\seg{AD}=b\,\sen{C}.
\]
Igualando las dos expresiones anteriores, se tiene:
\[
b\,\sen{C}=c\,\sen{B}.
\]
Dividiendo por $\sen{B}\sen{C}$, se tiene, finalmente:
\[
\frac{b}{\sen{B}}=\frac{c}{\sen{C}}
\]
Si se trazan perpendiculares a los lados $b$ y $c$ por los v'ertices $B$ y $C$, se demuestran las dem'as proporciones indicadas en la ecuaci'on \eqref{E1} empleando el mismo razonamiento expuesto anteriormente. De esta forma, el teorema queda demostrado.

\begin{figure}[htb]
\centering\fbox{\includegraphics[scale=0.8]{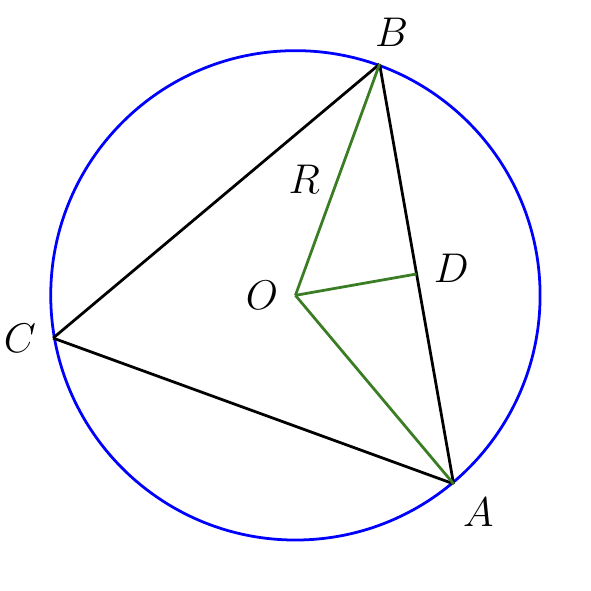}}
\caption{Construcci'on geom'etrica para encontrar la constante de proporcionalidad del teorema del seno.}\label{F3}
\end{figure}

Para terminar esta discusi'on, se buscar'a una interpretaci'on geom'etrica de la constante de proporcionalidad del teorema del seno. Para ello, la figura \ref{F3} ser'a de mucha utilidad.

Sea el tri'angulo $ABC$, en el cual se  circunscribe una circunferencia de radio $R$ y centro en $O$. Se trazan l'ineas desde $O$ hasta $A$ y $B$ y se traza una perpendicular por $O$ al lado $AB$, cort'andolo en $D$. Entonces:

El tri'angulo $OAB$ es is'osceles, porque $\seg{OB}=\seg{OA}=R$, y como $OD\bot AB$, entonces $\angle{BOD}=\angle{DOA}$.

Como el 'angulo $\angle{ACB}$ es inscrito, entonces es la mitad del 'angulo $\angle{AOB}$,  y como $\angle{BOD}=\angle{DOA}$, se tiene que $\angle{ACB}=\angle{BOD}=\angle{DOA}$. 

Del tri'angulo $OAD$, $\seg{AD}=\seg{AB}/2=R\sen{\angle{DOA}}=
R\sen{\angle{ACB}}$, y por lo tanto
\[
\frac{\seg{AB}}{\sen{\angle{ACB}}}=2R.
\]
Del teorema del seno, se deduce que
\begin{equation}\label{E2}
	\frac{\seg{AB}}{\sen{\angle{ACB}}}=
	\frac{\seg{BC}}{\sen{\angle{BAC}}}=
	\frac{\seg{CA}}{\sen{\angle{CBA}}}=2R
\end{equation}
La ecuaci'on \eqref{E2} permite concluir que \textit{la constante de proporcionalidad del teorema del seno aplicado a un tri'angulo oblicu'angulo es el di'ametro de la circunferencia circunscrita a 'el. }
\subsection{El teorema del coseno}\label{S2-2}
El teorema del coseno dice lo siguiente: {\it El cuadrado de un lado cualquiera en un tri'angulo oblicu'angulo es igual a la suma de los cuadrados de los lados restantes menos dos veces el producto de estos por el coseno del 'angulo comprendido entre ellos.}

Este teorema se escribe, de acuerdo al tri'angulo de la figura \ref{F1}, como:
\begin{align}
	a^2&=b^2+c^2-2bc\cos{A}\label{E3}\\
	b^2&=c^2+a^2-2ca\cos{B}\label{E4}\\
	c^2&=a^2+b^2-2ab\cos{C}\label{E5}
\end{align}
La demostraci'on de este teorema, al igual que la del teorema del seno, se hace en dos etapas:

Caso I: El tri'angulo es acut'angulo. De la figura \ref{F2}-a, al aplicar el teorema de Pit'agoras a los tri'angulos $ACD$ y $ABD$, se tiene:
\begin{align}
b^2&=\seg{AD}^2+\seg{DC}^2\label{E6}\\
c^2&=\seg{AD}^2+\seg{BD}^2\label{E7}
\end{align}
Restando la ecuaci'on \eqref{E7} de la ecuaci'on \eqref{E6} miembro a miembro:
\begin{equation}\label{E8}
	b^2-c^2=\seg{DC}^2-\seg{BD}^2
\end{equation}
Del tri'angulo $ABD$:
\[
\seg{BD}=c\cos{B},
\]
y como $\seg{DC}=a-\seg{BD}=a-c\cos{B}$, la ecuaci'on \eqref{E8} puede reescribirse de la siguiente forma:
\[
b^2-c^2=(a-c\cos{B})^2-c^2\cos^2{B} 
\]
Despejando $b^2$ de la ecuaci'on anterior y simplificando, se llega finalmente a:
\[
b^2=c^2+a^2-2ca\cos{B}
\]
Caso II: El tri'angulo es obtus'angulo. De la figura \ref{F2}-b, al aplicar el teorema de Pit'agoras a los tri'angulos $ACD$ y $ABD$, se tiene:
\begin{align}
b^2&=\seg{AD}^2+\seg{DC}^2\label{E9}\\
c^2&=\seg{AD}^2+\seg{BD}^2\label{E10}
\end{align}
Restando la ecuaci'on \eqref{E10} de la ecuaci'on \eqref{E9} miembro a miembro:
\begin{equation}\label{E11}
	b^2-c^2=\seg{DC}^2-\seg{BD}^2
\end{equation}
Del tri'angulo $ABD$:
\[
\seg{BD}=c\cos{\angle{ABD}}=c\cos(180^\circ-B)
\]
y como $\seg{DC}=a+\seg{BD}=a+c\cos{(180^\circ-B)}$, la ecuaci'on \eqref{E11} puede reescribirse de la siguiente forma:
\[
b^2-c^2=[a+c\cos{(180^\circ-B)}]^2-c^2\cos^2{(180^\circ-B)}=a^2+2ac\cos{(180^\circ-B)} 
\]
Despejando $b^2$ de la ecuaci'on anterior y teniendo en cuenta que $\cos{(180^\circ-B)}=-\cos{B}$, se llega finalmente a:
\[
b^2=c^2+a^2-2ca\cos{B}
\]
Las ecuaciones \eqref{E3} y \eqref{E5} se obtienen al aplicar el razonamiento anterior a los tri'angulos de la figura \ref{F2} si se les trazan perpendiculares a los lados $b$ y $c$ por los v'ertices $B$ y $C$. Por lo tanto, el teorema queda demostrado.
\subsection{El teorema de la tangente}\label{S2-3}
Este teorema, conocido tambi'en como las analog'ias de N'eper para tri'angulos planos, dice lo siguiente: {\it La suma de dos lados cualesquiera de un tri'angulo oblicu'angulo es a su diferencia como la tangente de la semisuma de sus 'angulos opuestos a estos lados es a la tangente de la semidiferencia de estos 'angulos.}

Matem'aticamente, el teorema se expresa, de acuerdo al tri'angulo de la figura \ref{F1}, como:
\begin{align}
	\frac{a+b}{a-b}&=\cfrac{\tg{\cfrac{A+B}{2}}}{\tg{\cfrac{A-B}{2}}}&
	\frac{b+a}{b-a}&=\cfrac{\tg{\cfrac{B+A}{2}}}{\tg{\cfrac{B-A}{2}}}\label{E12}\\
	\frac{b+c}{b-c}&=\cfrac{\tg{\cfrac{B+C}{2}}}{\tg{\cfrac{B-C}{2}}}&
	\frac{c+b}{c-b}&=\cfrac{\tg{\cfrac{C+B}{2}}}{\tg{\cfrac{C-B}{2}}}\label{E13}\\
	\frac{c+a}{c-a}&=\cfrac{\tg{\cfrac{C+A}{2}}}{\tg{\cfrac{C-A}{2}}}&
	\frac{a+c}{a-c}&=\cfrac{\tg{\cfrac{A+C}{2}}}{\tg{\cfrac{A-C}{2}}}\label{E14}
\end{align}
La demostraci'on de este teorema puede hacerse tanto anal'iticamente como geom'etricamente. 

{\it Demostraci'on anal'itica}: Del teorema del seno, se sabe que
\[
\frac{a}{\sen{A}}=\frac{b}{\sen{B}}
\]
Cambiando los medios:
\[
\frac{a}{b}=\frac{\sen{A}}{\sen{B}}
\]
Sumando y restando uno a ambos miembros de la ecuaci'on anterior:
\[
\frac{a}{b}+1=\frac{\sen{A}}{\sen{B}}+1\quad\quad\frac{a}{b}-1=\frac{\sen{A}}{\sen{B}}-1
\]
Dividiendo las ecuaciones anteriores miembro a miembro y simplificando, se tiene:
\[
	\frac{a+b}{a-b}=\frac{\sen{A}+\sen{B}}{\sen{A}-\sen{B}}
\]
Expresando la suma y diferencia de senos del miembro derecho de la ecuaci'on anterior como productos de senos y cosenos, se tiene\footnote{Esta operaci'on se hace mediante las f'ormulas de transformaci'on:
\begin{align*}
	\sen{x}+\sen{y}&=2\sen{\frac{x+y}{2}}\cos{\frac{x-y}{2}}\\
	\sen{x}-\sen{y}&=2\sen{\frac{x-y}{2}}\cos{\frac{x+y}{2}}
\end{align*}
}:
\[
	\frac{a+b}{a-b}=\cfrac{2\sen{\cfrac{A+B}{2}}\cos{\cfrac{A-B}{2}}}
	{2\sen{\cfrac{A-B}{2}}\cos{\cfrac{A+B}{2}}}
\]
Teniendo en cuenta que $\tg{x}=\dfrac{\sen{x}}{\cos{x}}$, se tiene finalmente:
\[
\frac{a+b}{a-b}=\cfrac{\tg{\cfrac{A+B}{2}}}{\tg{\cfrac{A-B}{2}}}
\]
La segunda relaci'on de la ecuaci'on \eqref{E12} se deduce de la ecuaci'on anterior teniendo en cuenta que $a-b=-(b-a)$ y $\tg{\dfrac{A-B}{2}}=-\tg{\dfrac{B-A}{2}}$, no sin antes multiplicar por $-1$ la ecuaci'on en cuesti'on. Las dem'as relaciones se obtienen al aplicar este procedimiento a las dem'as proporciones dadas por el teorema del seno.

En la mayor'ia de textos sobre trigonometr'ia, este teorema se obtiene de la forma expuesta anteriormente, lo cual no quiere decir que no pueda deducirse mediante m'etodos geom'etricos. La prueba que se presenta a continuaci'on\footnote{Esta demostraci'on se tom'o de \cite{B1}, p'ag. 43--44.}, tiene la ventaja de encontrar geom'etricamente unas ecuaciones que ser'an 'utiles m'as adelante. 
\begin{figure}[htb]
\centering
\fbox{\includegraphics[scale=0.8]{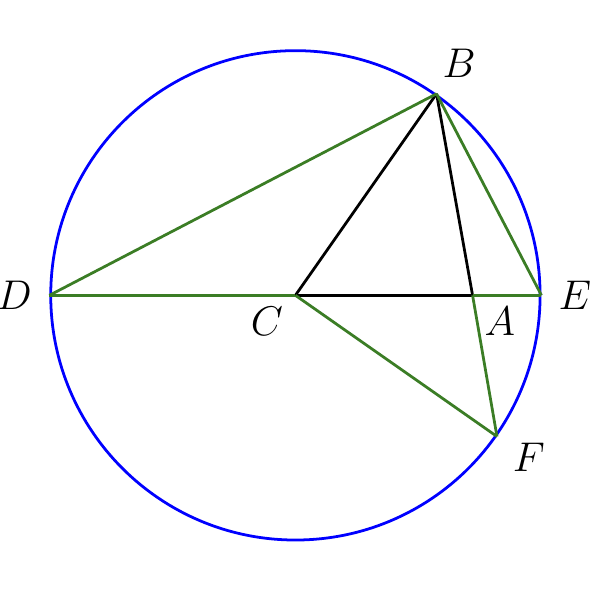}}
\caption{Construcci'on geom'etrica para la demostraci'on del teorema de la tangente.}\label{F4}
\end{figure}

{\it Demostraci'on geom'etrica}: Dado el tri'angulo $ABC$ (ver figura \ref{F4}), donde $a>b$, se traza una circunferencia de radio $\seg{BC}=a$ y centro en $C$. Luego, se prolonga la recta $AB$ por $A$ hasta que corte la circunferencia en $F$. Tambi'en se prolonga la recta $AC$ por ambos extremos hasta que alcance la circunferencia en $D$ y $E$. Finalmente, se une $C$ con $F$, $B$ con $D$, $E$ y $F$. 

De la figura en cuesti'on, se tiene que:

El 'angulo $\angle DBE$ es recto, pues est'a inscrito a la semicircunferencia $DFE$.

El 'angulo $\angle BCD$, al ser externo al tri'angulo $ABC$, es igual a la suma de los 'angulos $A$ y $B$; y como subtiende el arco $BD$, entonces el 'angulo $\angle BED$ inscrito a este arco es la mitad del 'angulo $\angle BCD$. En resumen:%\footnote{Aqu'i, se comete el abuso de notaci'on con respecto a la nomenclatura de los 'angulos internos del tri'angulo a tratar mencionado en la secci'on \ref{C1}. El lector no debe perder esto de vista.}
\begin{equation}\label{E15}
	A+B=\angle BCD=2\angle BED
\end{equation}
Por otro lado, al ser el 'angulo $A$ externo al tri'angulo $ACF$, es igual a la suma de los 'angulos $\angle ECF$ y $\angle AFC$; y como $\seg{CB}=\seg{CF}$, por ser radios de la circunferencia, el tri'angulo $BCF$ es is'osceles, lo que implica que $\angle AFC=B$ y por ende $\angle ECF=A-\angle AFC=A-B$. Adem'as, como el 'angulo $\angle ECF$ subtiende el arco $EF$, entonces el 'angulo $EBF$ inscrito a este arco es la mitad del 'angulo $\angle ECF$. De esto se concluye que
\begin{equation}\label{E16}
	A-B=\angle ECF=2\angle EBF
\end{equation}
De las ecuaciones \eqref{E15} y \eqref{E16}, se obtiene:
\begin{align}
	\angle BED&=\frac{A+B}{2}\label{E17}\\
	\angle EBF&=\frac{A-B}{2}\label{E18}
\end{align}
Los 'angulos $\angle ADB$ y $\angle ABD$ son complementarios de los 'angulos $\angle BED$ y $\angle EBF$, respectivamente.

Aplicando el teorema del seno al tri'angulo $ABD$:
\[
	\frac{\seg{AD}}{\sen{\angle ABD}}=\frac{\seg{AB}}{\sen{\angle ADB}}
\]
Como $\seg{AD}=\seg{AC}+\seg{CD}=\seg{AC}+\seg{BC}=a+b$, $\seg{AB}=c$, $\angle ABD=90-\angle EBF$ y $\angle ADB=90-\angle BED$, la ecuaci'on anterior se convierte en
\begin{equation}\label{E19}
\cfrac{a+b}{\cos{\cfrac{A-B}{2}}}=\cfrac{c}{\cos{\cfrac{A+B}{2}}}
\end{equation}
Aplicando el teorema del seno al tri'angulo $ABE$:
\[
	\frac{\seg{AE}}{\sen{\angle EBF}}=\frac{\seg{AB}}{\sen{\angle BED}}
\]
Como $\seg{AE}=\seg{CE}-\seg{AC}=\seg{BC}-\seg{AC}=a-b$ y $\seg{AB}=c$, la ecuaci'on anterior se convierte en
\begin{equation}\label{E20}
\cfrac{a-b}{\sen{\cfrac{A-B}{2}}}=\cfrac{c}{\sen{\cfrac{A+B}{2}}}
\end{equation}
Dividiendo las ecuaciones \eqref{E19} y \eqref{E20} miembro a miembro, se tiene:
\[
\cfrac{a+b}{\cos{\cfrac{A-B}{2}}}\cdot\cfrac{\sen{\cfrac{A-B}{2}}}{a-b}=
\cfrac{c}{\cos{\cfrac{A+B}{2}}}\cdot\cfrac{\sen{\cfrac{A+B}{2}}}{c}
\]
Despu'es de simplificar y agrupar t'erminos, se obtiene finalmente:
\[
\frac{a+b}{a-b}=\cfrac{\tg{\cfrac{A+B}{2}}}{\tg{\cfrac{A-B}{2}}}
\]
Con respecto a esta demostraci'on, se debe hacer las siguientes observaciones: 

1) Las relaciones deducidas a lo largo de la prueba son independientes de si el tri'angulo es acut'angulo u obtus'angulo. Sin embargo, la construcci'on depende de la condici'on $a>b$, lo cual implica que, para que sea completa, se debe construir un tri'angulo oblicu'angulo tal que $a<b$, y la construcci'on auxiliar comienza trazando una circunferencia de radio $b$ con centro en $C$. De ah'i en adelante, tanto la finalizaci'on de los trazos auxiliares como la demostraci'on son similares a la ya presentada. Si el lector se anima a completarla, debe llegar a lo siguiente:
\[
\frac{b+a}{b-a}=\cfrac{\tg{\cfrac{B+A}{2}}}{\tg{\cfrac{B-A}{2}}},
\]
Si bien, la ecuaci'on anterior y la deducida sean equivalentes desde el punto de vista algebraico, geom'etricamente no lo son. En dicho escenario, la primera solo es v'alida si $a>b$ y la segunda lo ser'a si $a<b$. Aunque esto sea un detalle insignificante en la pr'actica, no lo es en la teor'ia. %El lector debe tener siempre en mente este tipo de cosas.
  
2) Si se tiene en cuenta que $(A+B)/2=90-C/2$, las ecuaciones \eqref{E19} y \eqref{E20} se pueden reescribir de la siguiente manera:
\begin{align}
	\frac{a+b}{c}&=\cfrac{\cos{\cfrac{A-B}{2}}}{\sen{\cfrac{C}{2}}}
	&\frac{a-b}{c}&=\cfrac{\sen{\cfrac{A-B}{2}}}{\cos{\cfrac{C}{2}}}\label{E21}
\end{align}
Estas ecuaciones son un par del conjunto de relaciones conocidas como las f'ormulas de Mollweide, las cuales hacen parte de un conjunto de ecuaciones usadas para verificar que los par'ametros de un tri'angulo obtenidos a partir de otros sean los correctos. En la secci'on \ref{S2-6} se deducir'an estas ecuaciones anal'iticamente.

3) Las demostraciones presentadas dan la sensaci'on de que este teorema es una consecuencia directa o indirecta del teorema del seno. Sin embargo, esto no es as'i. En la literatura se encuentran demostraciones de este teorema que recurren a ingeniosas construcciones auxiliares en sus razonamientos, en las cuales no se necesita el teorema del seno. El lector interesado en algunas de ellas, puede consultar los siguientes textos: \cite{B2}, p'ag. 138--139; \cite{B3}, p'ag. 37 y  \cite{B4}, p'ag. 47--48.
\subsection{Las f'ormulas del 'angulo medio}\label{S2-4}
Si $2p=a+b+c$, entonces en todo tri'angulo oblicu'angulo:
\begin{align}
	\cos{\frac{A}{2}}=\sqrt{\frac{p(p-a)}{bc}}&&
	\sen{\frac{A}{2}}=\sqrt{\frac{(p-b)(p-c)}{bc}}\label{E22}\\
	\cos{\frac{B}{2}}=\sqrt{\frac{p(p-b)}{ac}}&&
	\sen{\frac{B}{2}}=\sqrt{\frac{(p-a)(p-c)}{ac}}\label{E23}\\
	\cos{\frac{C}{2}}=\sqrt{\frac{p(p-c)}{ab}}&&
	\sen{\frac{C}{2}}=\sqrt{\frac{(p-a)(p-b)}{ab}}\label{E24}
\end{align}
\textit{Demostraci'on:} Del teorema del coseno, ecuaci'on \eqref{E3}:
\[
\cos{A}=\frac{b^2+c^2-a^2}{2bc}
\]
Como $2\cos^2{\frac{A}{2}}=1+\cos{A}$ y $2\sen^2{\frac{A}{2}}=1-\cos{A}$, se tiene:
\begin{align*}
	\cos^2{\frac{A}{2}}&=\frac{1+\cos{A}}{2}=\cfrac{1+\cfrac{b^2+c^2-a^2}{2bc}}{2}=
	\frac{(a+b+c)(b+c-a)}{4bc}\\
	\sen^2{\frac{A}{2}}&=\frac{1-\cos{A}}{2}=\cfrac{1-\cfrac{b^2+c^2-a^2}{2bc}}{2}=
	\frac{(a+c-b)(a+b-c)}{4bc}\\
\end{align*}
Como $2p=a+b+c$, entonces $b+c-a=2(p-a)$, $a+c-b=2(p-b)$ y $a+b-c=2(p-c)$; las ecuaciones anteriores se convierten en:
\begin{align*}
	\cos{\frac{A}{2}}&=\sqrt{\frac{(a+b+c)(b+c-a)}{4bc}}=\sqrt{\frac{p(p-a)}{bc}}\\
	\sen{\frac{A}{2}}&=\sqrt{\frac{(a+c-b)(a+b-c)}{4bc}}=\sqrt{\frac{(p-b)(p-c)}{bc}}\\
\end{align*} 
En lo anterior, solo se tom'o el signo positivo de la raiz cuadrada porque $A/2$ siempre es agudo. Las dem'as relaciones se obtienen aplicando el procedimiento anterior a las restantes relaciones dadas del teorema del coseno.

\begin{figure}[htb!]
\centering\fbox{\includegraphics[scale=0.8]{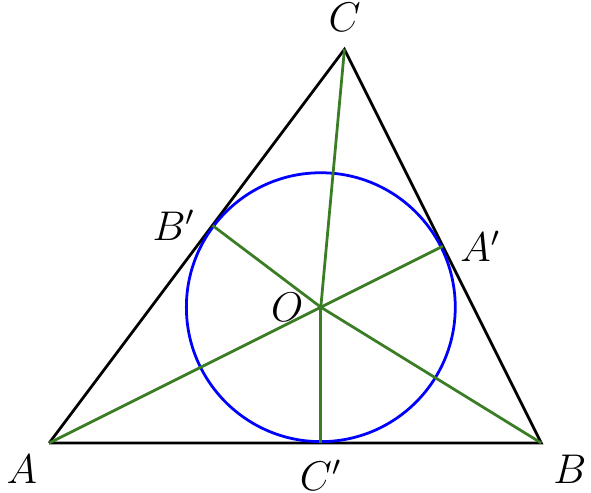}}
\caption{Construcci'on geom'etrica para la deducci'on de las f'ormulas del 'angulo medio}\label{F5}
\end{figure}

De las ecuaciones \eqref{E22}, \eqref{E23} y \eqref{E24}, se obtienen las tangentes del 'angulo medio:
\begin{align*}
	\tg{\frac{A}{2}}&=\sqrt{\frac{(p-b)(p-c)}{bc}}\sqrt{\frac{bc}{p(p-a)}}=
	\frac{1}{p-a}\sqrt{\frac{(p-a)(p-b)(p-c)}{p}}\\
	\tg{\frac{B}{2}}&=\sqrt{\frac{(p-a)(p-c)}{ac}}\sqrt{\frac{ac}{p(p-b)}}=
	\frac{1}{p-b}\sqrt{\frac{(p-a)(p-b)(p-c)}{p}}\\
	\tg{\frac{C}{2}}&=\sqrt{\frac{(p-a)(p-b)}{ab}}\sqrt{\frac{ab}{p(p-c)}}=
	\frac{1}{p-c}\sqrt{\frac{(p-a)(p-b)(p-c)}{p}}
\end{align*}
Si denotamos como $r$ a la expresi'on
\begin{equation}\label{E25}
	r=\sqrt{\frac{(p-a)(p-b)(p-c)}{p}},
\end{equation}
las ecuaciones anteriores pueden reescribirse de la siguiente manera:
\begin{align}
	\tg{\frac{A}{2}}&=\frac{r}{p-a}\label{E26}\\
	\tg{\frac{B}{2}}&=\frac{r}{p-b}\label{E27}\\
	\tg{\frac{C}{2}}&=\frac{r}{p-c}\label{E28}
\end{align}
En este 'ultimo conjunto de ecuaciones, $r$ es el radio de la circunferencia inscrita al tri'angulo $ABC$. Las f'ormulas del 'angulo medio, en los textos de trigonometr'ia, se demuestran anal'iticamente. Sin embargo, al igual que con el teorema de la tangente, pueden encontrarse vali'endose de la geometr'ia. Para ello, se recurrir'a a la figura \ref{F5}.

Al tri'angulo $ABC$, se trazan las bisectrices a los 'angulos $A$, $B$ y $C$ hasta que se encuentren en $O$, el incentro del tri'angulo en cuesti'on. Luego, se trazan perpendiculares a $AB$, $BC$ y $CA$ por $O$ cuyos pies son $C'$, $A'$ y $B'$ respectivamente. Finalmente, se traza la circunferencia inscrita al tri'angulo, cuyo radio se llamar'a $r$. De la figura se deduce:

$\seg{OA'}=\seg{OB'}=\seg{OC'}=r$, por ser radios de la circunferencia inscrita al tri'angulo $ABC$. 

Los tri'angulos $AB'O$ y $AC'O$ son rect'angulos en $B'$ y $C'$ respectivamente, pues por construcci'on, $OB'\perp AC$ y $OC'\perp AB$. Como $\seg{OA}=\seg{OA}$ y $\seg{OB'}=\seg{OC'}$, entonces, estos tri'angulos son iguales, y como consecuencia $\seg{AB'}=\seg{AC'}$.

Los tri'angulos $BC'O$ y $A'BO$ son rect'angulos en $C'$ y $A'$ respectivamente, pues por construcci'on, $OB'\perp AC$ y $OC'\perp AB$. Como $\seg{OB}=\seg{OB}$ y $\seg{OA'}=\seg{OC'}$, entonces, estos tri'angulos son iguales y $\seg{BC'}=\seg{BA'}$.

Los tri'angulos $B'CO$ y $A'CO$ son rect'angulos en $B'$ y $A'$ respectivamente, pues por construcci'on, $OB'\perp AC$ y $OA'\perp BC$. Como $\seg{OC}=\seg{OC}$ y $\seg{OA'}=\seg{OB'}$, entonces, estos tri'angulos son iguales y $\seg{CB'}=\seg{CA'}$.

Si el per'imetro se denota como $2p$, entonces:
\[
	\seg{AC'}+\seg{BC'}+\seg{BA'}+\seg{CA'}+\seg{CB'}+\seg{AB'}=2p,
\]
y teniendo en cuenta que $\seg{AB'}=\seg{AC'}$, $\seg{BC'}=\seg{BA'}$ y $\seg{CB'}=\seg{CA'}$, se tiene:
\begin{equation}\label{E29}
	\seg{AB'}+\seg{BC'}+\seg{CA'}=p.
\end{equation}
Por otro lado:
\begin{align}
	c=&\seg{AC'}+\seg{BC'}=\seg{AB'}+\seg{BC'}\label{E30}\\
	b=&\seg{AB'}+\seg{CB'}=\seg{AB'}+\seg{CA'}\label{E31}\\
	a=&\seg{BA'}+\seg{CA'}=\seg{BC'}+\seg{CA'}\label{E32}
\end{align}
Combinando las ecuaciones \eqref{E29}, \eqref{E30}, \eqref{E31} y \eqref{E32}, se obtiene:
\begin{align}
	\seg{AB'}=\seg{AC'}=&p-a\label{E33}\\
	\seg{BC'}=\seg{BA'}=&p-b\label{E34}\\
	\seg{CB'}=\seg{CA'}=&p-c\label{E35}
\end{align}
De los tri'angulos $AB'O$, $BC'O$ y $CA'O$, se deduce:
\begin{align*}
	\tg{\frac{A}{2}}&=\frac{\,\seg{OB'}\,}{\,\seg{AB'}\,}=\frac{r}{p-a}\\
	\tg{\frac{B}{2}}&=\frac{\,\seg{OC'}\,}{\,\seg{BC'}\,}=\frac{r}{p-b}\\
	\tg{\frac{C}{2}}&=\frac{\,\seg{OA'}\,}{\,\seg{CA'}\,}=\frac{r}{p-c}
\end{align*}
Las ecuaciones anteriores son las tangentes de los 'angulos medios del tri'angulo $ABC$. Para encontrar el radio de la circunferencia inscrita, se aplica el teorema del seno al tri'angulo $ABO$:
\begin{equation}\label{E36}
	\seg{OA}=\cfrac{c\sen{\cfrac{B}{2}}}{\sen{\GP{180^\circ-\cfrac{A+B}{2}}}}=
	\cfrac{c\sen{\cfrac{B}{2}}}{\cos{\cfrac{C}{2}}},
\end{equation}
y tambi'en se aplica el teorema del seno al tri'angulo $ACO$:
\begin{equation}\label{E37}
	\seg{OA}=\cfrac{b\sen{\cfrac{C}{2}}}{\sen{\GP{180^\circ-\cfrac{A+C}{2}}}}=
	\cfrac{c\sen{\cfrac{C}{2}}}{\cos{\cfrac{B}{2}}}
\end{equation}
Multiplicando las ecuaciones \eqref{E36} y \eqref{E37}, y empleando las relaciones encontradas para las tangentes de los 'angulos medios, se tiene:
\begin{equation}\label{E38}
	\seg{OA}^2=bc\tg{\frac{B}{2}}\tg{\frac{C}{2}}=\frac{bcr^2}{(p-b)(p-c)}
\end{equation}
Aplicando el teorema de Pit'agoras al tri'angulo $AB'O$:
\[
\seg{OA}^2=\seg{AB'}^2+\seg{OB'}^2\quad\therefore\quad\frac{bcr^2}{(p-b)(p-c)}=(p-a)^2+r^2
\] 
Despejando $r^2$:
\[
	r^2=(p-a)^2\,{\frac{(p-b)(p-c)}{bc-(p-b)(p-c)}}=
	(p-a)^2\,{\frac{(p-b)(p-c)}{p(b+c-p)}}=
	(p-a)^2\,{\frac{(p-b)(p-c)}{p(p-a)}}.
\]
De la ecuaci'on anterior se encuentra que $r$ es igual a:
\[
r=\sqrt{{\frac{(p-a)(p-b)(p-c)}{p}}},
\]
Expresi'on que coincide con la ecuaci'on \eqref{E25}. Teniendo en cuenta las ecuaciones \eqref{E25}, \eqref{E33} y \eqref{E38}, se deduce del tri'angulo $AB'O$:
\begin{align*}
	\cos{\frac{A}{2}}=&\frac{\,\seg{AB'}\,}{\,\seg{OA}\,}=(p-a)\sqrt{\frac{(p-b)(p-c)}{bcr^2}}
	=\sqrt{\frac{p(p-a)}{bc}}\\
	\sen{\frac{A}{2}}=&\frac{\,\seg{OB'}\,}{\,\seg{OA}\,}=r\sqrt{\frac{(p-b)(p-c)}{bcr^2}}
	=\sqrt{\frac{(p-b)(p-c)}{bc}}
\end{align*}
Las relaciones para los 'angulos restantes se encuentran examinando los tri'angulos $BC'O$ y $CA'O$, no sin antes de haber encontrado relaciones similares a $\seg{OA}$ para $\seg{OB}$ y $\seg{OC}$.

La demostraci'on anterior se bas'o sobre un tri'angulo acut'angulo. Para un tri'angulo obtus'angulo, tanto la construcci'on como la cadena de razonamientos usados no cambia, por ello, no se presenta aqu'i.
\subsection{Las f'ormulas de la proyecci'on}\label{S2-5}
En todo tri'angulo oblicu'angulo:
\begin{align}
a&=b\cos{C}+c\cos{B}\label{E39}\\
b&=c\cos{A}+a\cos{C}\label{E40}\\
c&=a\cos{B}+b\cos{A}\label{E41}
\end{align}
Demostraci'on: La prueba se divide, como de costumbre, en dos casos:

Caso I: El tri'angulo es acut'angulo. De la figura \ref{F2}-a: $a=\seg{BD}+\seg{DC}$. De los tri'angulos $ABD$ y $ACD$, se tiene que $\seg{BD}=c\cos{B}$ y $\seg{DC}=b\cos{C}$, respectivamente. De lo anterior se obtiene que
\[
a=\seg{BD}+\seg{DC}=b\cos{C}+c\cos{B}
\]
Caso II: El tri'angulo es obtus'angulo. De la figura \ref{F2}-b: $a=\seg{DC}-\seg{DB}$. De los tri'angulos $ABD$ y $ACD$, se tiene que $\seg{DB}=c\cos{(180-B)}$ y $\seg{DC}=b\cos{C}$, respectivamente. De lo anterior se obtiene que
\[
a=\seg{DC}-\seg{DB}=b\cos{C}-c\cos{(180-B)}=b\cos{C}+c\cos{B}
\] 
Si se trazan perpendiculares a los lados $b$ y $c$ por los v'ertices $B$ y $C$, y siguiendo el razonamiento anterior, se llegan a las dem'as relaciones. Gracias a ello, el teorema queda demostrado.
\subsection{Las f'ormulas de Mollweide}\label{S2-6}
En todo tri'angulo oblicu'angulo:
\begin{align}
\frac{a+b}{c}&=\cfrac{\cos{{\cfrac{A-B}{2}}}}{\sen{\cfrac{C}{2}}}&
\frac{a-b}{c}&=\cfrac{\sen{{\cfrac{A-B}{2}}}}{\cos{\cfrac{C}{2}}}\label{E42}\\
\frac{b+a}{c}&=\cfrac{\cos{{\cfrac{B-A}{2}}}}{\sen{\cfrac{C}{2}}}&
\frac{b-a}{c}&=\cfrac{\sen{{\cfrac{B-A}{2}}}}{\cos{\cfrac{C}{2}}}\label{E43}\\
\frac{b+c}{a}&=\cfrac{\cos{{\cfrac{B-C}{2}}}}{\sen{\cfrac{A}{2}}}&
\frac{b-c}{a}&=\cfrac{\sen{{\cfrac{B-C}{2}}}}{\cos{\cfrac{A}{2}}}\label{E44}\\
\frac{c+b}{a}&=\cfrac{\cos{{\cfrac{C-B}{2}}}}{\sen{\cfrac{A}{2}}}&
\frac{c-b}{a}&=\cfrac{\sen{{\cfrac{C-B}{2}}}}{\cos{\cfrac{A}{2}}}\label{E45}\\
\frac{a+c}{b}&=\cfrac{\cos{{\cfrac{A-C}{2}}}}{\sen{\cfrac{B}{2}}}&
\frac{a-c}{b}&=\cfrac{\sen{{\cfrac{A-C}{2}}}}{\cos{\cfrac{B}{2}}}\label{E46}\\
\frac{c+a}{b}&=\cfrac{\cos{{\cfrac{C-A}{2}}}}{\sen{\cfrac{B}{2}}}&
\frac{c-a}{b}&=\cfrac{\sen{{\cfrac{C-A}{2}}}}{\cos{\cfrac{B}{2}}}\label{E47}
\end{align}
Demostraci'on: Al reescribir el teorema del seno de la siguiente manera:
\[
\frac{a}{c}=\frac{\sen{A}}{\sen{C}}\quad\quad\frac{b}{c}=\frac{\sen{B}}{\sen{C}},
\]
entonces, sumando y restando ambas ecuaciones miembro a miembro, se tiene:
\[
\frac{a+b}{c}=\frac{\sen{A}+\sen{B}}{\sen{C}}\quad\quad\frac{a-b}{c}=\frac{\sen{A}-\sen{B}}{\sen{C}}
\]
Transformando las sumas y diferencias de senos en productos y expresando $\sen{C}$ en t'erminos de $C/2$:
\[
\frac{a+b}{c}=\cfrac{2\sen{\cfrac{A+B}{2}}\cos{\cfrac{A-B}{2}}}{2\sen{\cfrac{C}{2}}\cos{\cfrac{C}{2}}}\quad\quad\frac{a-b}{c}=\cfrac{2\sen{\cfrac{A-B}{2}}\cos{\cfrac{A+B}{2}}}{2\sen{\cfrac{C}{2}}\cos{\cfrac{C}{2}}}
\]
Como $A+B=180-C$, entonces $\sen{\frac{A+B}{2}}=\cos{\frac{C}{2}}$ y $\cos{\frac{A+B}{2}}=\sen{\frac{C}{2}}$. Con esto en mente, se tiene finalmente que
\[
\frac{a+b}{c}=\cfrac{\cos{{\cfrac{A-B}{2}}}}{\sen{\cfrac{C}{2}}}\quad\quad
\frac{a-b}{c}=\cfrac{\sen{{\cfrac{A-B}{2}}}}{\cos{\cfrac{C}{2}}}
\]
De esta manera, se obtuvieron las ecuaciones \eqref{E42}. La primera de las ecuaciones \eqref{E43} se obtiene al cambiar el orden de los sumandos del numerador de la fracci'on del primer miembro y reemplazando $\cos{\dfrac{A-B}{2}}$ por $\cos{\dfrac{B-A}{2}}$ en la primera de las ecuaciones \eqref{E42}; mientras que la segunda surge al reemplazar respectivamente a $a-b$ y $\sen{\dfrac{A-B}{2}}$ por $-(b-a)$ y $-\sen{\dfrac{B-A}{2}}$ en la segunda de las ecuaciones \eqref{E43}, ecuaci'on a la cual se multiplico por $-1$ previamente. Las dem'as relaciones se encuentran mediante manipulaciones similares a las diversas proporciones dadas por el teorema del seno. De esta forma, el teorema queda demostrado. En la secci'on \ref{S2-3}, se indic'o c'omo se encontraban las ecuaciones \eqref{E42} de forma geom'etrica. Las ecuaciones \eqref{E43} se deducen para el caso $b>a$, el cual se dej'o como ejercicio para el lector.
\section{Soluci\'on de los diferentes casos}\label{S3}
Una vez establecidos los teoremas necesarios, ahora se mostrar'a su uso para resolver tri'angulos oblicu'angulos. La discusi'on de cada caso comienza con un tratamiento anal'itico detallado para el c'alculo de los par'ametros desconocidos y en metodolog'ias para corrobarlos. Al final de cada uno, se ilustra su soluci'on mediante construcciones geom'etricas.
\subsection{Primer caso: cuando se conocen dos 'angulos y un lado com'un a ellos}\label{S3-1}
Sup'ongase conocidos los 'angulos $B$ y $C$ junto con el lado $a$. Como la suma de los tres 'angulos internos es igual a dos rectos, el 'angulo restante ser'a el suplemento de la suma de los 'angulos conocidos, es decir:
\begin{equation}\label{E48}
A=180^{\circ}-(B+C).
\end{equation}
Una vez conocido el 'angulo $A$, los lados restantes se obtienen mediante el teorema del seno:
\begin{align}
\frac{b}{\sen{B}}&=\frac{a}{\sen{A}}\quad\therefore\quad b=a\,\frac{\sen{B}}{\sen{A}}\label{E49}\\
\frac{c}{\sen{C}}&=\frac{a}{\sen{A}}\quad\therefore\quad c=a\,\frac{\sen{C}}{\sen{A}}\label{E50}
\end{align}
Las ecuaciones \eqref{E48}, \eqref{E49} y \eqref{E50} proporcionan un solo valor para $A$, $b$ y $c$. Por lo tanto, el primer caso est'a resuelto. Para comprobar que los par'ametros est'an correctamente calculados, se puede acudir ya sea a las f'ormulas de la proyecci'on o las de Mollweide. Si se usan las primeras, la ecuaci'on \eqref{E39} ser'a la indicada; pero si se quieren emplear las segundas, las ecuaciones \eqref{E44} o \eqref{E45} son las id'oneas, siempre y cuando se reescriban de la siguiente manera:
\begin{align}
(b+c){\sen{\cfrac{A}{2}}}&={a}\,{\cos{{\cfrac{B-C}{2}}}}&
(b-c){\cos{\cfrac{A}{2}}}&={a}\,{\sen{{\cfrac{B-C}{2}}}}\label{E51}\\
(c+b){\sen{\cfrac{A}{2}}}&={a}\,{\cos{{\cfrac{C-B}{2}}}}&
(c-b){\cos{\cfrac{A}{2}}}&={a}\,{\sen{{\cfrac{C-B}{2}}}}\label{E52}
\end{align}
Cualquiera de las cuatro ecuaciones anteriores sirve para verificar la rectitud de los par'ametros calculados. El lector debe notar que la elecci'on de estas relaciones de Mollweide es totalmente arbitraria, lo cual permite escoger a su gusto cualquiera de las enunciadas en la secci'on donde fueron presentadas.

\begin{figure}[htb!]
\centering\fbox{\includegraphics[scale=0.8]{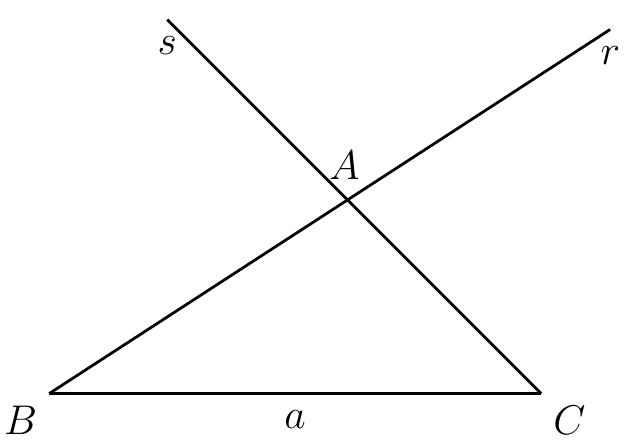}}
\caption{Construcci'on geom'etrica del tri'angulo para el primer caso.}\label{F6}
\end{figure}

Geom'etricamente, el tri'angulo se construye de la siguiente manera: se traza un segmento $\seg{BC}$ de longitud $a$. Por los extremos $B$ y $C$, se construyen semirrectas $r$ y $s$ tales que formen 'angulos $B$ y $C$ con el segmento trazado. Finalmente, se prolongan las semirectas hasta que se corten en un punto $A$, lo cual forma el tri'angulo $ABC$. Al cumplir este las condiciones del problema, se tiene finalmente el tri'angulo deseado (ver figura \ref{F6}).
\subsection{Segundo caso: cuando se conocen dos 'angulos y un lado opuesto a uno de ellos}\label{S3-2}
Sup'ongase conocidos los 'angulos $A$ y $B$ junto al lado $a$. De la propiedad de los 'angulos internos de un tri'angulo se deduce que el 'angulo restante es igual a:
\begin{equation}\label{E53}
C=180^{\circ}-(A+B).
\end{equation}
Los lados que faltan se calculan mediante la aplicaci'on del teorema del seno:
\begin{align}
\frac{b}{\sen{B}}&=\frac{a}{\sen{A}}\quad\therefore\quad b=a\,\frac{\sen{B}}{\sen{A}}\label{E54}\\
\frac{c}{\sen{C}}&=\frac{a}{\sen{A}}\quad\therefore\quad c=a\,\frac{\sen{C}}{\sen{A}}\label{E55}
\end{align}
Las ecuaciones \eqref{E53}, \eqref{E54} y \eqref{E55} proporcionan un solo valor para $C$, $b$ y $c$. Con esto se soluciona el segundo caso. Al igual que en el caso anterior, es indispensable corroborar los resultados obtenidos, lo cual se realiza mediante las f'ormulas de la proyecci'on (de las cuales, la ecuaci'on \eqref{E39} es la indicada) o usando una de las relaciones de Mollweide, ya sea las que est'an en las ecuaciones \eqref{E51} y \eqref{E52}, o cualquiera de las que est'an enlistadas en la secci'on \ref{S2-6}

\begin{figure}[htb!]
\centering\fbox{\includegraphics[scale=0.8]{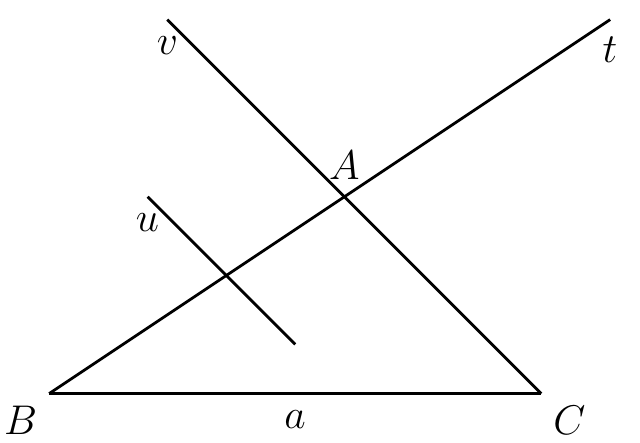}}
\caption{Construcci'on geom'etrica del tri'angulo para el segundo caso.}\label{F7}
\end{figure}

Para construir geom'etricamente el tri'angulo bajo estas condiciones, se emplea el siguiente proceder: se traza un segmento $\seg{BC}$ de longitud $a$. Por el extremo $B$, se construye la semirrecta $t$ tal que forme un 'angulo $B$ con el segmento trazado. En un punto cualquiera de $t$, se traza una recta $u$ tal que forme un 'angulo $A$ con ella. Hecho esto, se traza una recta $v$ paralela a $u$ que pase por $C$ y que cortar'a a $t$ en el punto $A$, dando origen al  tri'angulo $ABC$. Al cumplir este las condiciones del problema, se tiene finalmente el tri'angulo deseado (ver figura \ref{F7}). 
\subsection{Tercer caso: cuando se conocen dos lados y un 'angulo opuesto a uno de ellos}\label{S3-3}
Sup'ongase conocidos los lados $a$ y $b$ con el 'angulo $A$. Entonces el 'angulo opuesto al lado $B$ se obtiene a partir del teorema del seno:
\begin{equation}\label{E56}
\frac{a}{\sen{A}}=\frac{b}{\sen{B}}\quad\therefore\quad\sen{B}=\frac{b\sen{A}}{a}
\end{equation}
Dependiendo de los valores de $a$, $b$ y $A$, se puede tener uno, dos o ningun valor para $B$. Para ello, se analizan los siguientes casos:
\begin{enumerate}[i)]
\item  $a>b$: Sin importar si $A$ sea agudo u obtuso, de la ecuaci'on \eqref{E56} se deduce que $\sen{B}<1$, porque $a>b>b\sen{A}$, lo que significa que $B$ puede tomar dos valores, los cuales son suplementarios entre ellos. Como $a>b$, se tiene de las propiedades de los tri'angulos que $A>B$. De esto se deduce que si $A$ es agudo, $B$ tambi'en lo ser'a y si $A$ es obtuso, $B$ tiene que ser forzosamente agudo. Por lo tanto, se tendr'a un solo tri'angulo que cumpla las condiciones dadas.% \emph{Si $a>b$, entonces $B$ siempre ser'a agudo y su seno vendr'a dado por la ecuaci'on \eqref{E52}.}
\item $a<b$: En este caso, la ecuaci'on \eqref{E56} indica que $\sen{B}$ puede ser mayor, igual o menor que uno. En el primer escenario, $b\sen{A}>a$, lo cual implica que no existe un tri'angulo ; en el segundo escenario, si $A$ es agudo, $b\sen{A}=a$, y $B$ ser'a de $90^{\circ}$: se tendr'a un tri'angulo recto, pero si $A$ es obtuso, el problema no tiene soluci'on, pues al ser $A>B$, implica a su vez que $a>b$, contradicciendo as'i la restricci'on original; en el 'ultimo escenario, $b\sen{A}<a$, y $B$ puede tomar dos valores: como $a<b$  exige que $A<B$, entonces si $A$ es agudo se tendr'an dos tri'angulos cuyos 'angulos opuestos al lado $b$ son suplementarios y si $A$ es obtuso, se tendr'a un tri'angulo cuyo 'angulo opuesto al lado $b$ sea agudo.

\item $a=b$: En este caso, se tendr'a un tri'angulo si $A$ es agudo, el cual es is'osceles; y si $A$ es obtuso, el problema no tiene sentido, pues un tri'angulo no puede tener dos 'angulos obtusos.
\end{enumerate}%Un resumen de la discusi'on anterior se da en la tabla 1.
Una vez se determine $B$, el resto de par'ametros se encuentran de la siguiente forma:

Si solo se tiene un valor para $B$, entonces el 'angulo restante se deriva de la relaci'on entre los 'angulos internos de un tri'angulo:
\begin{equation}\label{E57}
C=180^{\circ}-(A+B);
\end{equation}
y el lado restante se obtiene del teorema del seno:
\begin{equation}\label{E58}
\frac{c}{\sen{C}}=\frac{a}{\sen{A}}\quad\therefore\quad c=a\,\frac{\sen{C}}{\sen{A}};
\end{equation}
Si se tienen dos valores para $B$, los cuales se denominar'an $B_1$ y $B_2$, entonces, los 'angulos restantes para cada tri'angulo se deducen de una forma similar a la ya presentada, es decir:
\begin{align}
C_1=180^{\circ}-(A+B_1),\label{E59}\\
C_2=180^{\circ}-(A+B_2);\label{E60}
\end{align}
y los lados restantes se obtienen del teorema del seno:
\begin{align}
c_1=a\,\frac{\sen{C_1}}{\sen{A}}\label{E61}\\
c_2=a\,\frac{\sen{C_2}}{\sen{A}}\label{E62}
\end{align}
Para verificar que los resultados obtenidos son los correctos, se usa una de las f'ormulas de Mollweide. Se debe insistir que su elecci'on es libre.

La construcci'on geom'etrica de este caso permite visualizar f'acilmente los casos en los cuales se tiene una, dos o ninguna soluci'on. Independientemente de si el 'angulo opuesto dado sea agudo u obtuso, la construcci'on, suponiendo dados $A$, $a$ y $b$, comienza al trazar el 'angulo $A$. En uno de sus lados, se ubica un punto $C$ tal que su distancia a su v'ertice sea $b$. Desde este punto, se baja una perpendicular al lado restante cuyo pie es $D$. Con centro en $C$ y radio $a$, se traza un arco tal que corte a la recta $AE$, intersecci'on que depender'a de la magnitud de $A$ y de la relaci'on de $a$ y $CD$.
\begin{figure}[htb!]
\centering
\begin{tabular}{cc}\fbox{
\includegraphics[scale=0.8]{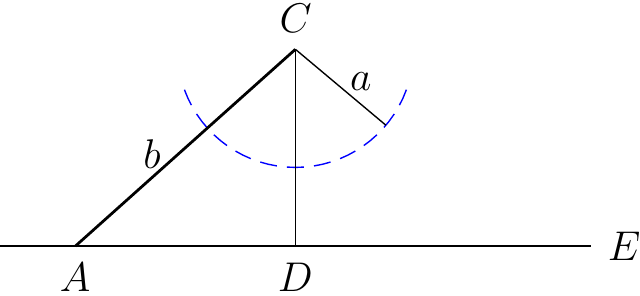}}
&\fbox{
\includegraphics[scale=0.8]{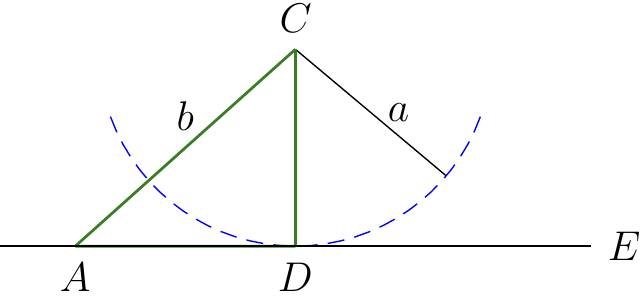}}
\\
a) & b)\\\fbox{
\includegraphics[scale=0.8]{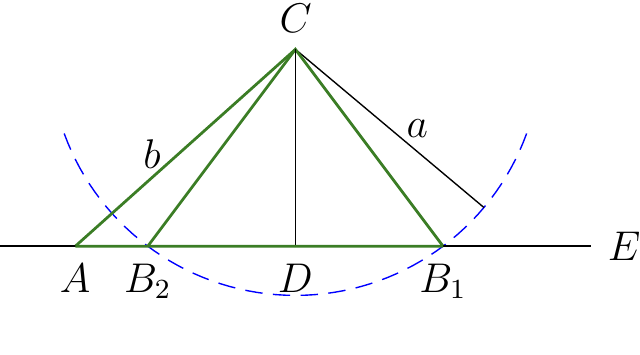}}
&\fbox{
\includegraphics[scale=0.8]{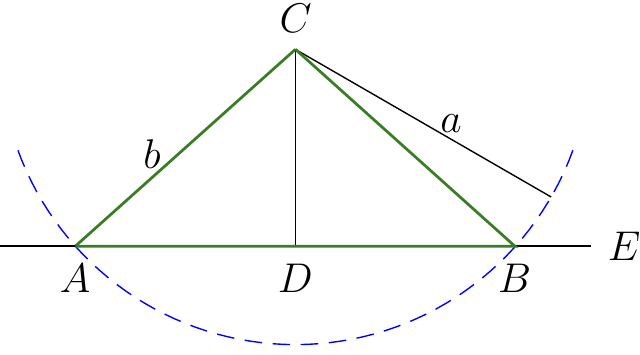}}
\\
c) & d)\\
\multicolumn{2}{c}{\fbox{
\includegraphics[scale=0.8]{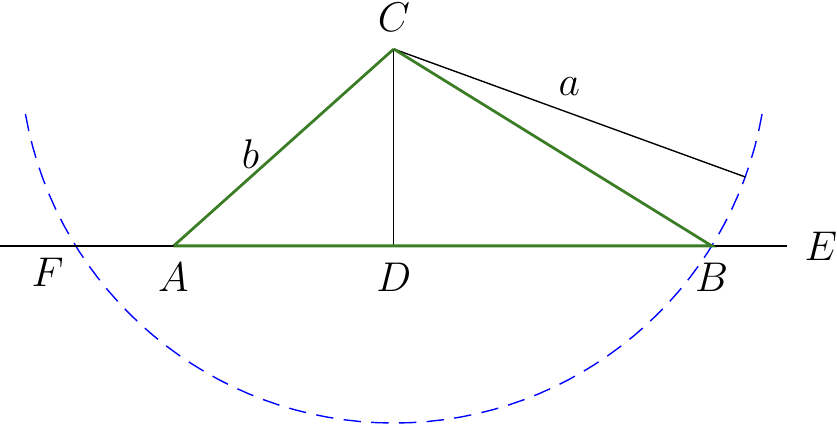}}
}\\
\multicolumn{2}{c}{e)}
\end{tabular}
\caption{Construcci\'on geom\'etrica del tercer caso cuando el \'angulo opuesto dado es agudo.}\label{F8}
\end{figure}

Cuando $A$ es agudo:

Si $a<\seg{CD}$, el arco no cortar'a a $AE$ y no se tendr'a soluci'on alguna (ver figura \ref{F8}-a);

Si $a=\seg{CD}$, el arco ser'a tangente a $AE$, y la soluci'on ser'a el tri'angulo $ACD$ (ver figura \ref{F8}-b) ;

Si $a>\seg{CD}$ y $a<b$, el arco cortar'a a $AE$ en los puntos $B_1$ y $B_2$. En este caso, se tendr'an dos soluciones materializadas en los tri'angulos $AB_1C$ y $AB_2C$ (ver figura \ref{F8}-c);

Si $a>\seg{CD}$ y $a=b$, el arco cortar'a a $AE$ en los puntos $A$ y $B$. En este caso, la soluci'on ser'a el tri'angulo is'osceles $ABC$ (ver figura \ref{F8}-d);

Si $a>\seg{CD}$ y $a>b$, el arco cortar'a a $AE$ en los puntos $B$ y $F$. Sin embargo, de los dos tri'angulos que podr'ian formarse, el tri'angulo  $ABC$ cumple con las condiciones impuestas por el segundo caso, lo cual lo convierte en la soluci'on buscada (ver figura \ref{F8}-e).
\begin{figure}[htb!]
\centering
\begin{tabular}{cc}
\centering\fbox{
\includegraphics[scale=0.8]{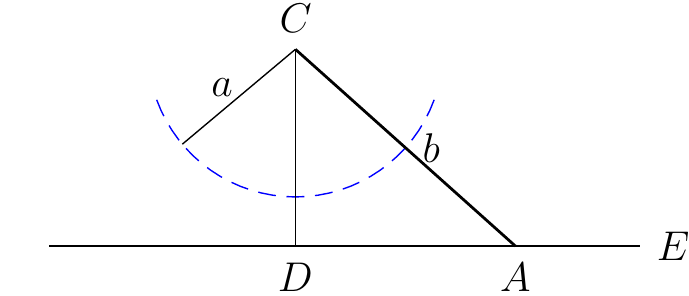}}
&\fbox{
\includegraphics[scale=0.8]{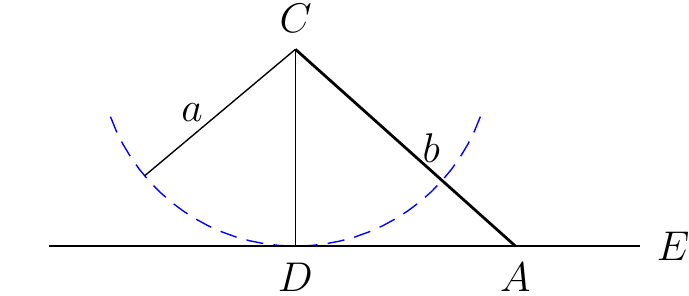}}
\\
a) & b)\\\fbox{
\includegraphics[scale=0.8]{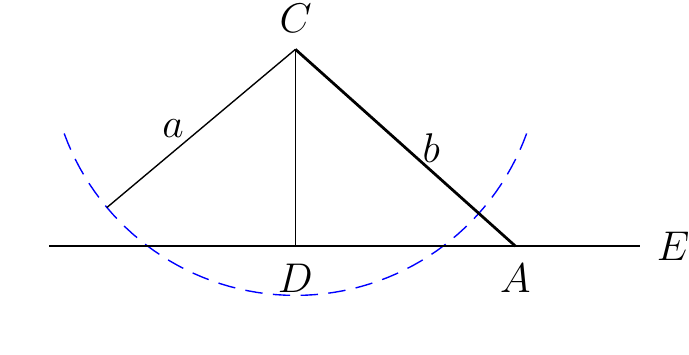}}
&\fbox{
\includegraphics[scale=0.8]{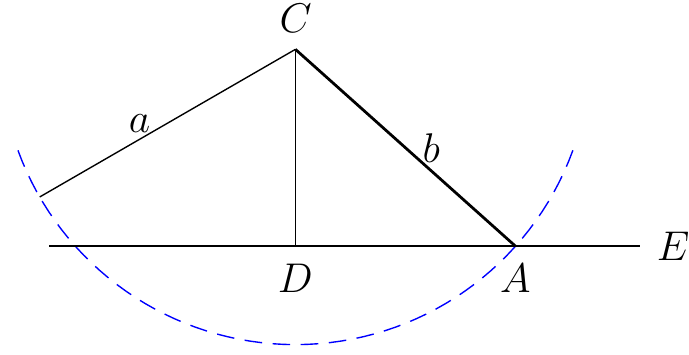}}
\\
c) & d)\\
\multicolumn{2}{c}{\fbox{
\includegraphics[scale=0.8]{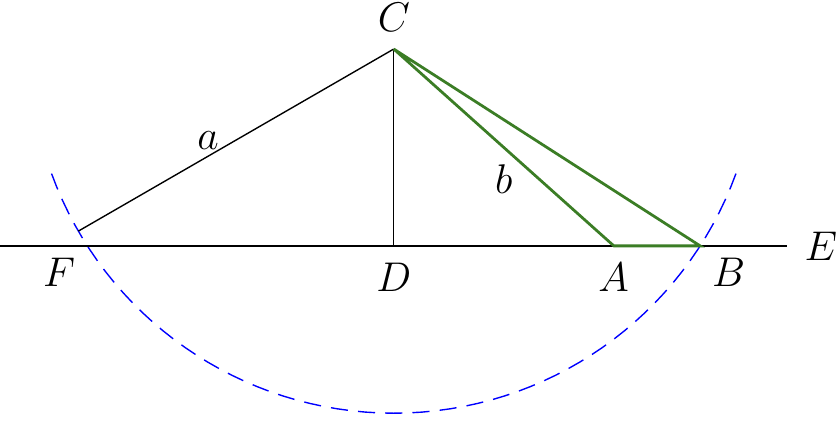}}
}\\
\multicolumn{2}{c}{e)}
\end{tabular}
\caption{Construcci'on geom'etrica del tercer caso cuando el 'angulo dado es obtuso.}\label{F9}
\end{figure}

Cuando $A$ es obtuso:

Ya sea que si $a<\seg{CD}$, $a=\seg{CD}$, $\seg{CD}<a<b$; y $\seg{CD}<a=b$, no existe un tri'angulo que cumpla las condiciones del caso tres, tal y como lo ilustran los cuadros a), b), c) y d) de la figura \ref{F9}. Sin embargo, cuando $a>b$, el arco en cuesti'on intersecta a la recta $AE$ en los puntos $B$ y $F$, de los cuales solo el primero hace parte del v'ertice del tri'angulo $ABC$. Este tri'angulo satisface los requisitos del tercer caso.

Con esta discusi'on geom'etrica, se corroboran las conclusiones obtenidas por medios anal'iticos en cuanto a la cantidad de soluciones que pueden presentarse de acuerdo a los valores de los par'ametros dados.
\subsection{Cuarto caso: cuando se conocen dos casos y el 'angulo comprendido entre ellos}\label{S3-4}
Sup'ongase dados los lados $a$ y $b$ con el 'angulo $C$. Entonces el lado $c$ se calcula usando el teorema del coseno:
\begin{equation}\label{E63}
c^{2}=a^{2}+b^{2}-2ab\cos{C}\quad\therefore\quad c=\sqrt{a^{2}+b^{2}-2ab\cos{C}}
\end{equation}
Los 'angulos restantes se calculan a partir del teorema del seno:
\begin{align}
\frac{a}{\sen{A}}&=\frac{c}{\sen{C}}\quad\therefore\quad\sen{A}=\frac{a\sen{C}}{c}\label{E64}\\
\frac{b}{\sen{B}}&=\frac{c}{\sen{C}}\quad\therefore\quad\sen{B}=\frac{b\sen{C}}{c}\label{E65}
\end{align}
Las ecuaciones \eqref{E64} y \eqref{E65} indican que los senos de $A$ y $B$ pueden ser mayores, iguales o menores que uno. Para el caso en cuesti'on, estos son menores que uno. Para asegurar esto, se puede ver en la figura \ref{F2}-a que, en el tri'angulo recto $ABD$, la hipotenusa $c$ es mayor que $\seg{BD}=b\sen{C}$, por ser el lado que opone el mayor 'angulo, lo cual garantiza que $\seg{BD}/c=\sen{B}$ sea menor que uno, si $C$ es 'agudo. Un razonamiento similar en la figura \ref{F2}-b revela que $\seg{BD}/c=\sen{B}$ es menor que uno para $C$ obtuso. Si en dichos tri'angulos se traza la altura con respecto a $B$, y efectuando un an'alisis paralelo al anterior, se concluir'a que $\sen{A}<1$, ya sea que $C$ sea agudo u obtuso.

Por otro lado, $\sen{A}$ y $\sen{B}$ proporcionan dos valores tanto para $A$ como para $B$, lo cual conducir'ia a m'as de una soluci'on. A continuaci'on, se muestra que siempre es posible escoger un valor para $A$ y uno para $B$. Si $C<90^{\circ}$, de la relaci'on de los 'angulos internos de un tri'angulo se deduce que $A+B>90^{\circ}$. Esta 'ultima desigualdad implica  tres opciones: $A$ y $B$ son agudos; $A$ es obtuso y $B$ es agudo o, $A$ es agudo y $B$ es obtuso. Para escoger la correcta, basta utilizar el hecho de que el mayor (menor) lado opone mayor (menor) 'angulo. Entonces, ordenando los lados de forma descendente (o ascendente), se elige los valores de los 'angulos opuestos a estos que cumplan dicho orden. Si $C>90^{\circ}$, entonces $A+B<90^{\circ}$, lo cual implica que $A$ y $B$ deben ser agudos. De esta forma, se asegura que el tercer caso tiene una soluci'on. Para comprobar que la soluci'on es correcta, la suma de los 'angulos internos debe ser de $180^{\circ}$.

El procedimiento expuesto anteriormente es quiza, el m'as empleado. Sin embargo, existe m'as de una forma de lidiar con este caso. El m'etodo expuesto a continuaci'on usa el teorema de la tangente. Un primer vistazo al teorema muestra que no es posible usarlo, pues no se conoce los 'angulos opuestos a los lados dados; pero, como se conoce el 'angulo comprendido entre ellos, de la relaci'on de los 'angulos internos de un tri'angulo, se deduce f'acilmente que el suplemento del 'angulo conocido es igual a la suma de los 'angulos restantes. En conclusi'on, el teorema de la tangente proporciona la semidiferencia de los 'angulos opuestos de los lados dados. Si se suponen conocidos los lados $a$, $b$ y el 'angulo $C$, entonces: 

La suma de los 'angulos iguales es igual a: $A+B=180^{\circ}-C$, y del teorema de la tangente, se tiene:
\begin{equation}\label{E66}
\frac{a+b}{a-b}=\cfrac{\tg{\cfrac{A+B}{2}}}{\tg{\cfrac{A-B}{2}}}
\quad\therefore\quad
\tg{\frac{A-B}{2}}=\frac{a-b}{a+b}\cotg{\frac{C}{2}}
\end{equation}
La semidiferencia de los 'angulos opuestos ser'a positiva y estar'a comprendida entre $0$ y $90^{\circ}$ si $a>b$; ser'a negativa si $a<b$ y estar'a comprendida entre $-90^{\circ}$ y $0$; y ser'a nula si $a=b$. Si se denomina a $\delta$ a la semidiferencia en cuesti'on, los 'angulos se obtienen as'i:
\begin{align}
A=&\frac{A+B}{2}+\frac{A-B}{2}=\frac{180^{\circ}-C}{2}+\delta\label{E67}\\
B=&\frac{A+B}{2}-\frac{A-B}{2}=\frac{180^{\circ}-C}{2}-\delta\label{E68}
\end{align}
El lado $c$ se obtiene a partir del teorema del seno:
\begin{equation}\label{E69} 
c=a\,\frac{\sen{C}}{\sen{A}}=b\,\frac{\sen{C}}{\sen{B}}
\end{equation}
Para verificar que los resultados son correctos para este m'etodo expuesto, se puede usar una de las relaciones de Mollweide. 

\begin{figure}[htb!]
\centering\fbox{\includegraphics[scale=0.8]{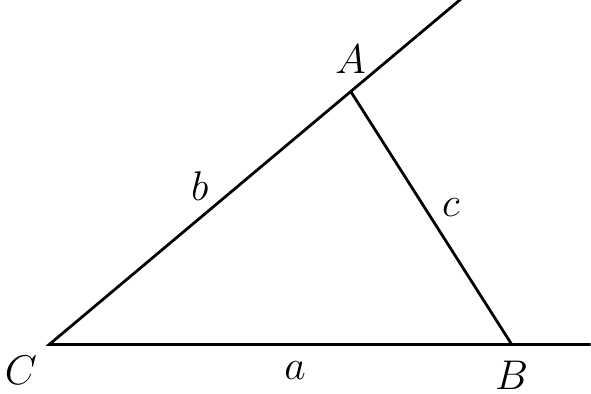}}
\caption{Construcci'on geom'etrica del tri'angulo para el cuarto caso.}\label{F10}
\end{figure}

Para finalizar, se muestra c'omo se construye el tri'angulo conocidos dos lados $a$ y $b$ y el 'angulo comprendido entre ellos $C$: se traza el 'angulo $C$, y a partir de su v'ertice, se ubican sobres sus lados dos puntos $A$ y $B$ tales que $\seg{AC}=b$ y $\seg{BC}=a$. Hecho esto, se une con una recta los puntos en cuesti'on. El tri'angulo $ABC$ es el deseado (ver figura \ref{F10}).
\subsection{Quinto caso: cuando se conocen los tres lados}\label{S3-5}
Si se conocen los tres lados, entonces del teorema del coseno:
\begin{align}
	a^2&=b^2+c^2-2bc\cos{A}\quad\therefore\quad\cos{A}=\frac{b^2+c^2-a^2}{2bc}\label{E70}\\
	b^2&=c^2+a^2-2ca\cos{B}\quad\therefore\quad\cos{B}=\frac{c^2+a^2-b^2}{2ca}\label{E71}\\
	c^2&=a^2+b^2-2ab\cos{C}\quad\therefore\quad\cos{C}=\frac{a^2+b^2-c^2}{2ab}\label{E72}
\end{align}
Para mirar cu'ando se tiene una soluci'on o ninguna, se analizar'a la ecuaci'on \eqref{E70}. Si existe una soluci'on, entonces $\cos{A}$ es menor que uno pero mayor a menos uno:
\begin{equation}\label{E73} 
-1<{\frac{b^2+c^2-a^2}{2bc}}<1
\end{equation}
Despu'es de algunas manipulaciones, la desigualdad anterior se puede escribir como:
\[
(b-c)^{2}<a^{2}<(b+c)^{2}
\]
Tomando ra'iz cuadrada a los miembros de la desigualdad anterior, se tiene:
\begin{equation}\label{E74} 
|b-c|<a<b+c
\end{equation}
La relaci'on anterior plasma la conocida relaci'on entre los lados de un tri'angulo: La suma de dos lados cualesquiera de un tri'angulo es mayor que el tercero, pero menor que su diferencia. Un an'alisis similar para $\cos{B}$ y $\cos{C}$ conduce a la misma conclusi'on. En resumen, dados los tres lados de un tri'angulo, se tendr'a una soluci'on si se satisfacen las desigualdades:
\begin{align}
|b-c|<&a<b+c	\label{E75}\\
|c-a|<&b<c+a	\label{E76}\\
|a-b|<&c<a+b	\label{E77}
\end{align}
Otra forma de resolver este caso es usando las f'ormulas del 'angulo medio. Si solo se desea conocer los 'angulos, entonces:

Se calcula el semiper'imetro $p=(a+b+c)/2$ y luego, se puede escoger las f'ormulas que proporcionan los senos de los 'angulos medios:
\begin{align}
	\sen{\frac{A}{2}}=&\sqrt{\frac{(p-b)(p-c)}{bc}}\label{E78}\\
	\sen{\frac{B}{2}}=&\sqrt{\frac{(p-a)(p-c)}{ac}}\label{E79}\\
	\sen{\frac{C}{2}}=&\sqrt{\frac{(p-a)(p-b)}{ab}}\label{E80},
\end{align}
o las que dan los cosenos de los 'angulos medios:
\begin{align}
	\cos{\frac{A}{2}}=&\sqrt{\frac{p(p-a)}{bc}}\label{E81}\\
	\cos{\frac{B}{2}}=&\sqrt{\frac{p(p-b)}{ac}}\label{E82}\\
	\cos{\frac{C}{2}}=&\sqrt{\frac{p(p-c)}{ab}}\label{E83}.
\end{align}

Si adem'as de los 'angulos, se quiere saber el radio de la circunferencia inscrita al tri'angulo que se quiere resolver, entonces se usan las ecuaciones de las tangentes de los 'angulos medios\footnote{En los textos que empleaban logaritmos para hacer los c'alculos, las f'ormulas de las tangentes eran las predilectas para resolver este caso, pues ellas est'an en funci'on de $p$, $p-a$, $p-b$ y $p-c$. Esto implica que solo era necesario encontrar los logaritmos de estas cantidades y con ellos, el c'alculo de $r$ y de los 'angulos internos era f'acil.}:
\begin{align}
	\tg{\frac{A}{2}}&=\frac{r}{p-a}\label{E84}\\
	\tg{\frac{B}{2}}&=\frac{r}{p-b}\label{E85}\\
	\tg{\frac{C}{2}}&=\frac{r}{p-c}\label{E86}
\end{align}
donde $r$ est'a dado por la ecuaci'on \eqref{E25}. Si se quiere saber si se puede resolver o no el tri'angulo, las ecuaciones anteriores indican que el semiper'imetro debe ser mayor que los tres lados. Esta condici'on es m'as sencilla que la dada por las relaciones \eqref{E75}, \eqref{E76} y\eqref{E77}. Para verificar que los resultados son correctos para los dos m'etodos expuestos, la suma de los tres 'angulos debe ser de $180^{\circ}$.

\begin{figure}[htb!]
\centering
\begin{tabular}{ccc}\fbox{\includegraphics[scale=0.8]{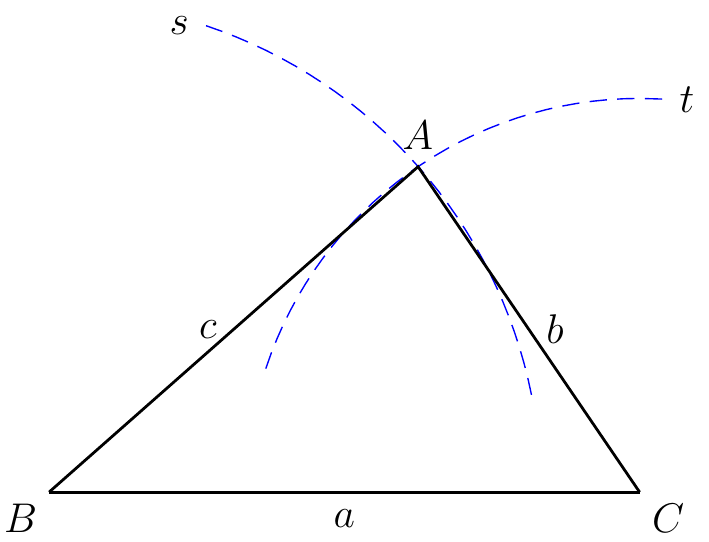}}
&&\fbox{\includegraphics[scale=0.8]{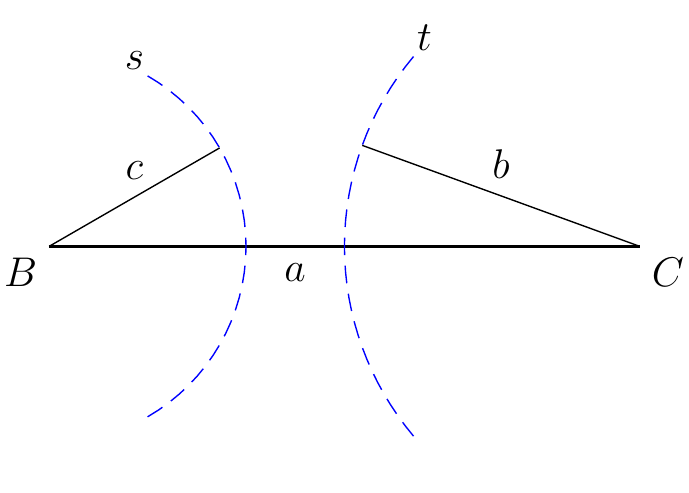}}
\\
a) && b)
\end{tabular}
\caption{Construcci'on geom'etrica del tri'angulo para el quinto caso: a) cuando la suma de los dos lados es mayor que el tercero; b) cuando la suma de los dos lados es menor o igual al lado restante.}\label{F11}
\end{figure}

La construcci'on geom'etrica para este caso es la siguiente: Se traza un segmento $\seg{BC}$ de longitud $a$. Se trazan arcos $s$ y $t$ hacia un mismo lado de $BC$ con centros en $B$ y $C$ y radios respectivos $c$ y $b$. Si estos arcos se intersectan en un punto $A$, el tri'angulo formado al unir los puntos $A$, $B$ y $C$ es el buscado, lo cual suceder'a si la suma de los radios de los arcos es mayor que la separaci'on de sus centros (ver figura \ref{F11}-a). Si esto no llegase a suceder, no se obtendr'a un tri'angulo que tenga por lados de longitud $a$, $b$ y $c$ (ver figura \ref{F11}-b).
\subsection{Otros casos de resoluci'on de tri'angulos}\label{S3-6}
Adem'as de los ya expuestos, existen otros cuyos par'ametros dados pueden ser la suma de dos lados, uno o dos 'angulos exteriores, el per'imetro o el 'area del tri'angulo, las alturas, las bisectrices, etc. A continuaci'on, se analizar'an dos casos en particular: dado los 'angulos internos y el per'imetro; y dado un lado, su correspondiente 'angulo opuesto y la suma de los lados restantes.
\subsubsection{Dado los 'angulos internos y el per'imetro:} Este caso puede resolverse de dos maneras. La primera consiste en usar las ecuaciones de las tangentes de los 'angulos medios. Utilizando las ecuaciones \eqref{E25}--\eqref{E28}, se forman las siguientes relaciones:
\begin{align}
\tg{\frac{A}{2}}\tg{\frac{B}{2}}&=\frac{r^2}{(p-a)(p-b)}=1-\frac{c}{p}\label{E87}\\
\tg{\frac{B}{2}}\tg{\frac{C}{2}}&=\frac{r^2}{(p-b)(p-c)}=1-\frac{a}{p}\label{E88}\\
\tg{\frac{C}{2}}\tg{\frac{A}{2}}&=\frac{r^2}{(p-c)(p-a)}=1-\frac{b}{p}\label{E89}
\end{align}
Despejando $a$, $b$, y $c$ de las ecuaciones anteriores, se tiene:
\begin{align}
a&=p\GP{1-\tg{\frac{B}{2}}\tg{\frac{C}{2}}}\label{E90}\\
b&=p\GP{1-\tg{\frac{A}{2}}\tg{\frac{C}{2}}}\label{E91}\\
c&=p\GP{1-\tg{\frac{A}{2}}\tg{\frac{B}{2}}}\label{E92}
\end{align}
Para probar que los resultados dados por las ecuaciones \eqref{E90}, \eqref{E91} y \eqref{E92} son correctos, basta con comprobar que $a+b+c=2p$. La segunda forma de atacar este problema es usando el teorema del seno. Empleando las propiedades de las proporciones\footnote{Si $\frac{a}{b}=\frac{c}{d}=\frac{e}{f}=\cdots$, entonces
\[
\frac{a}{b}=\frac{c}{d}=\frac{e}{f}=\cdots=\frac{a+c+e+\cdots}{b+d+f+\cdots}
\]}, se tiene que
\begin{equation}\label{E93}
\frac{a}{\sen{A}}=\frac{b}{\sen{B}}=\frac{c}{\sen{C}}=
\frac{a+b+c}{\sen{A}+\sen{B}+\sen{C}}=\frac{2p}{\sen{A}+\sen{B}+\sen{C}}
\end{equation}
De este conjunto de ecuaciones se obtiene:
\begin{align}
a&=\frac{2p\sen{A}}{\sen{A}+\sen{B}+\sen{C}}\label{E94}\\
b&=\frac{2p\sen{B}}{\sen{A}+\sen{B}+\sen{C}}\label{E95}\\
c&=\frac{2p\sen{C}}{\sen{A}+\sen{B}+\sen{C}}\label{E96}
\end{align}
Si se utiliza la identidad trigonom'etrica
\begin{equation}\label{E97}
\sen{A}+\sen{B}+\sen{C}=4\cos{\frac{A}{2}}\cos{\frac{B}{2}}\cos{\frac{C}{2}},\quad A+B+C=180^\circ
\end{equation}
Las ecuaciones \eqref{E94}, \eqref{E95} y \eqref{E96} se pueden rescribir de la siguiente manera:
\begin{align}
a&=\frac{2p\cdot2\sen{\dfrac{A}{2}}\cos{\dfrac{A}{2}}}{4\cos{\dfrac{A}{2}}\cos{\dfrac{B}{2}}\cos{\dfrac{C}{2}}}=\frac{p\sen{\dfrac{A}{2}}}{\cos{\dfrac{B}{2}}\cos{\dfrac{C}{2}}}\label{E98}\\
b&=\frac{2p\cdot2\sen{\dfrac{B}{2}}\cos{\dfrac{B}{2}}}{4\cos{\dfrac{A}{2}}\cos{\dfrac{B}{2}}\cos{\dfrac{C}{2}}}=\frac{p\sen{\dfrac{B}{2}}}{\cos{\dfrac{A}{2}}\cos{\dfrac{C}{2}}}\label{E99}\\
c&=\frac{2p\cdot2\sen{\dfrac{C}{2}}\cos{\dfrac{C}{2}}}{4\cos{\dfrac{A}{2}}\cos{\dfrac{B}{2}}\cos{\dfrac{C}{2}}}=\frac{p\sen{\dfrac{C}{2}}}{\cos{\dfrac{A}{2}}\cos{\dfrac{B}{2}}}\label{E100}
\end{align}
El resultado es correcto si $a+b+c=2p$. El teorema del seno reduce este problema a uno de reparto proporcional.
\subsubsection{Dado un lado, su correspondiente 'angulo opuesto y la suma de los lados restantes:}
Sup'ongase conocidos $a$, $A$ y $b+c$. De la formula de Mollweide
\[
\frac{b+c}{a}=\cfrac{\cos{\cfrac{B-C}{2}}}{\sen{\cfrac{A}{2}}}
\]
se deduce que\footnote{El signo $\pm$ en la ecuaci'on \eqref{E101} refleja que se desconoce si $B$ es mayor o menor que $C$.}
\begin{equation}\label{E101}
\frac{B-C}{2}=\pm\arccos{\GP{\frac{b+c}{a}\sen{\frac{A}{2}}}}
\end{equation}
Como $\dfrac{B+C}{2}=90^\circ-\dfrac{A}{2}$, al combinar esta relaci'on con la ecuaci'on \eqref{E101} se obtiene:
\begin{align}
B&=90^\circ-\dfrac{A}{2}\pm\arccos{\GP{\frac{b+c}{a}\sen{\frac{A}{2}}}}\label{E102}\\
C&=90^\circ-\dfrac{A}{2}\mp\arccos{\GP{\frac{b+c}{a}\sen{\frac{A}{2}}}}\label{E103}
\end{align}
Estas ecuaciones indican que se tienen dos soluciones: un tri'angulo cuyos 'angulos adyacentes al lado dado son 
\begin{equation}\label{E104}
B_1=90^\circ-\dfrac{A}{2}+\arccos{\GP{\frac{b+c}{a}\sen{\frac{A}{2}}}}\quad
C_1=90^\circ-\dfrac{A}{2}-\arccos{\GP{\frac{b+c}{a}\sen{\frac{A}{2}}}}
\end{equation}
y el otro tiene como 'angulos adyacentes al susodicho lado
\begin{equation}\label{E105}
B_2=90^\circ-\dfrac{A}{2}-\arccos{\GP{\frac{b+c}{a}\sen{\frac{A}{2}}}}\quad
C_2=90^\circ-\dfrac{A}{2}+\arccos{\GP{\frac{b+c}{a}\sen{\frac{A}{2}}}}
\end{equation}
Al comparar las ecuaciones \eqref{E104} y \eqref{E105}, se tiene que $B_1=C_2$ y $B_2=C_1$. Como ambos tri'angulos tienen un lado respectivo igual y sus correspondientes 'angulos adyacentes a ese lado iguales, por congruencia de tri'angulos, estos son iguales. Esto implica que se tiene una sola soluci'on. Por lo tanto, los 'angulos buscados son:
\begin{align}
B&=90^\circ-\dfrac{A}{2}+\arccos{\GP{\frac{b+c}{a}\sen{\frac{A}{2}}}}\label{E106}\\
C&=90^\circ-\dfrac{A}{2}-\arccos{\GP{\frac{b+c}{a}\sen{\frac{A}{2}}}}\label{E107}
\end{align}
N'otese que su elecci'on es arbitraria, pudiendose elegir la otra soluci'on. Una vez determinado los 'angulos, sus lados opuestos respectivos se obtienen mediante el teorema del seno:
\begin{align}
b&=\frac{a\sen{B}}{\sen{A}}\label{E108}\\
c&=\frac{a\sen{C}}{\sen{A}}\label{E109}
\end{align}
Las ecuaciones \eqref{E106} y \eqref{E107} muestran que no siempre se tiene una soluci'on para todos los valores de $a$, $b+c$ y $A$. Primero, debe tenerse en cuenta que siempre $b+c>a$. Esto significa que el termino $\dfrac{b+c}{a}\sen{\dfrac{A}{2}}$ no siempre es menor o igual uno. Luego, para tener una soluci'on, se debe cumplir que
\begin{equation}\label{E110}
\frac{b+c}{a}\sen{\frac{A}{2}}\leq1\quad\therefore\quad\sen{\frac{A}{2}}\leq\frac{a}{b+c}
\end{equation}
La construcci'on geom'etrica de este caso descansa en que, independientemente del valor de $A$  que satisfaga la desigualdad anterior y al ser la suma de los lados desconocidos constante, todos los v'ertices opuestos al lado dado est'an en una elipse cuyos focos son los extremos de dicho lado. Dicha elipse tiene como semieje mayor $a'=\dfrac{b+c}{2}$, semidistancia focal $c'=\dfrac{a}{2}$ y semieje menor $b'=\sqrt{a'^2-c'^2}$. Entonces, si se logra encontrar geom'etricamente un punto de la elipse tal que sus radio vectores formen un 'angulo $A$, el tri'angulo podr'a construirse.

\begin{figure}[htb!]
\centering\fbox{
\includegraphics[scale=0.8]{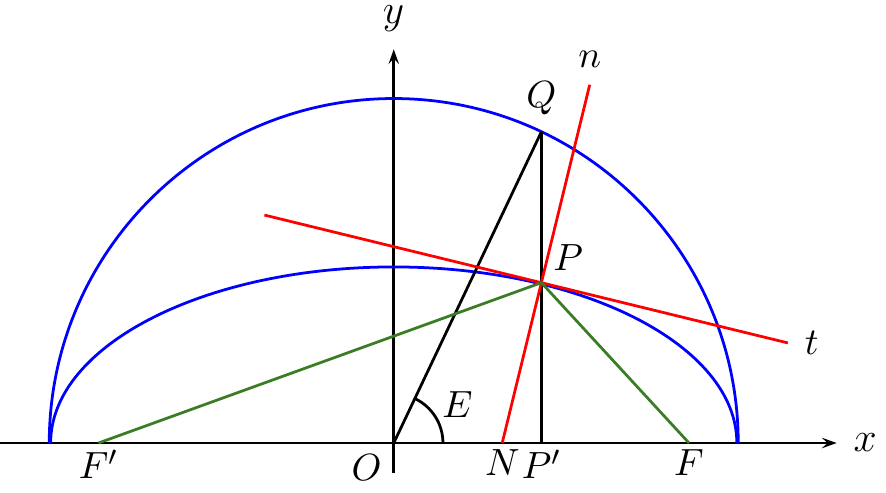}}
\caption{Construcci'on para determinar una relaci'on entre el 'angulo $A$ y la anomal'ia exc'entrica de la elipse asociada al tri'angulo buscado}\label{F12}
\end{figure}

Para ello, se establecer'a una relaci'on entre este 'angulo y un par'ametro de la elipse que sea f'acil de interpretar, en este caso, su anomal'ia exc'entrica. Con este fin, se ubica la elipse en un sistema cartesiano de coordenadas tal que el origen est'e en el punto medio del lado $BC$ y los ejes $x$  y $y$ coincidir'a y ser'a perpendicular  a 'el, respectivamente. Adicional a ello, se traza su circunferencia auxiliar. Lo descrito anteriormente se muestra en la figura \ref{F12}.

\begin{figure}[htb!]
\centering
\begin{tabular}{cc}\fbox{
\includegraphics[scale=0.725]{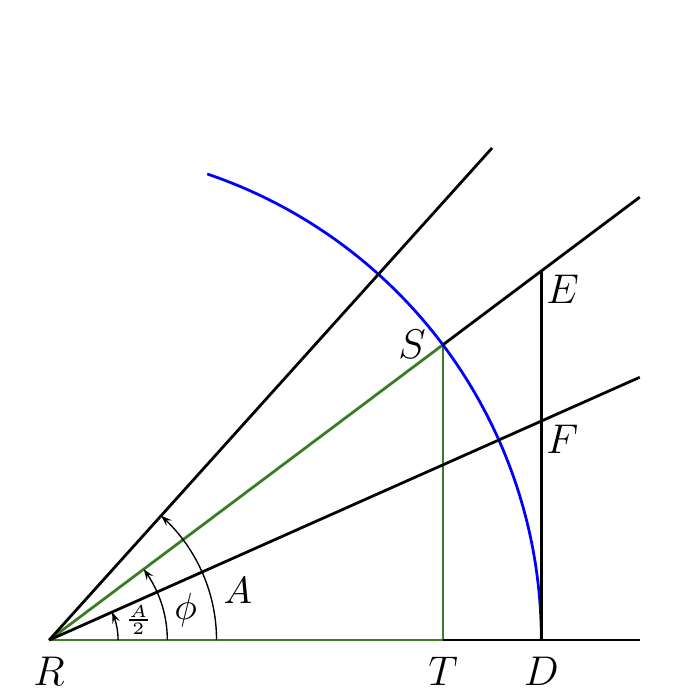}}&\fbox{
\includegraphics[scale=0.725]{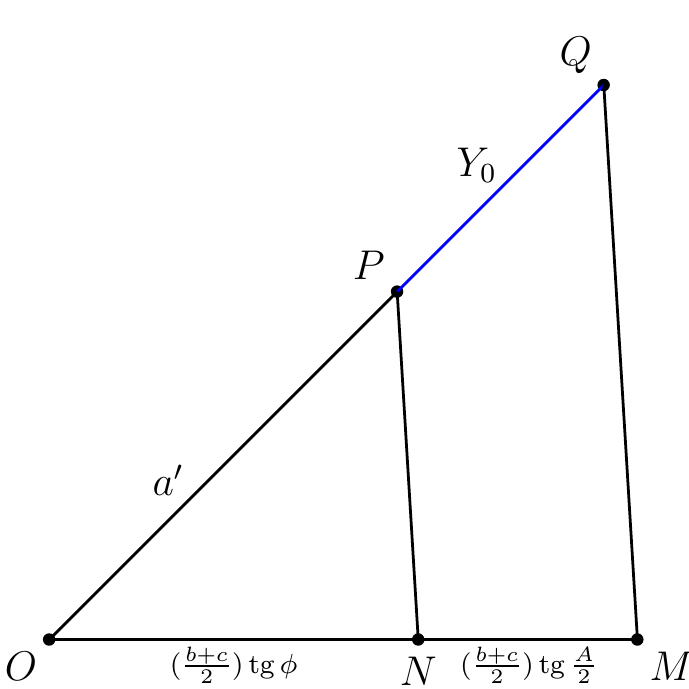}}\\
a) & b)\\
\multicolumn{2}{c}{\fbox{\includegraphics[scale=0.725]{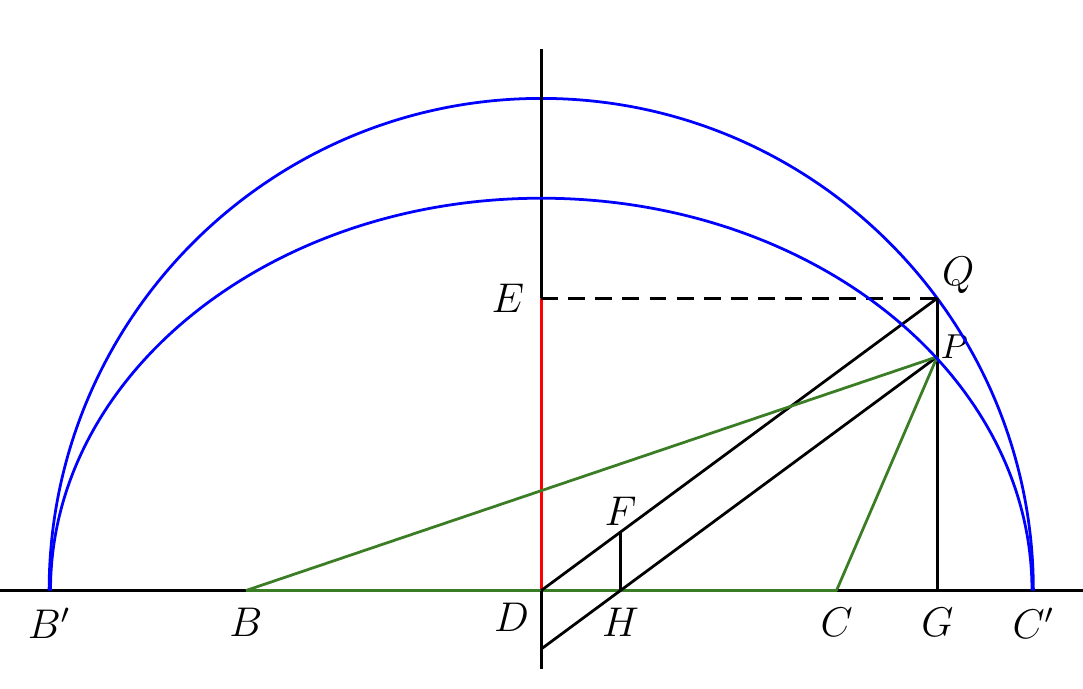}}}\\
\multicolumn{2}{c}{c)}
\end{tabular}
\caption{Construcci'on geom'etrica del tri'angulo dados un lado, su 'angulo opuesto y la suma de los lados restantes: a) Determinaci'on de los segmentos de longitud $(\frac{b+c}{2})\tg{\phi}$ y $(\frac{b+c}{2})\tg{\frac{\al}{2}}$; b) Construcci'on de la ordenada del punto $Q$; c) Ubicaci'on del v'ertice restante.}\label{F13}
\end{figure}

Sea $P(x_0,y_0)$ un punto de la elipse y $Q$ su punto asociado a la circunferencia auxiliar. Sea $F$ y $F'$ los focos de la elipse cuyas coordenadas son $(c',0)$ y $(-c',0)$; $n$ y $t$ son respectivamente, la recta normal y tangente a la elipse en el punto $P$. El 'angulo $\angle FPF'$ se denotar'a como $\al$ y el 'angulo $\angle QON$ es la anomal'ia exc'entrica $E$ correspondiente al punto $P$. Se sabe por geometr'ia anal'itica que las pendientes de las rectas normal $m_n$ y tangente $m_t$ a la elipse en $P$ son\footnote{Aunque estos resultados pueden obtenerse mediante las t'ecnicas del c'alculo diferencial, tambi'en pueden derivarse usando m'etodos algebraicos. En \cite{B5}, se expone un m'etodo basado en usar ecuaciones param'etricas de la recta para analizar su intersecci'on con una c'onica. En particular, en el cap'itulo 8, se efectua dicho an'alisis para la elipse.}:
\begin{align}
m_t&=-\frac{b'^2x_0}{a'^2y_0}\label{E111}\\
m_n&=-\frac{1}{m_t}=\frac{a'^2y_0}{b'^2x_0}\label{E112}
\end{align}
Las pendientes de las rectas $PF$ y $PF'$ son:
\begin{align}
m_{PF}&=\cfrac{y_0}{x_0-\cfrac{a}{2}}\label{E113}\\
 m_{PF'}&=\cfrac{y_0}{x_0+\cfrac{a}{2}}\label{E114}
\end{align}
Se sabe por las propiedades de la elipse que el 'angulo que forma los radio vectores de una elipse en un punto de ella lo biseca la recta normal a la elipse en ese punto\footnote{Vease, \cite{B7}, p'ag. 187, \cite{B6}, p'ag. 182.}. Con ello, se puede calcular la tangente de la mitad de ese 'angulo a partir de las pendientes de la normal y de uno de los radio vectores, es decir:
\[
\tg{\frac{\al}{2}}=\frac{m_n-m_{PF'}}{1+m_nm_{PF'}}=
\cfrac{\cfrac{a'^2y_0}{b'^2x_0}-\cfrac{y_0}{x_0+\cfrac{a}{2}}}{1+\cfrac{a'^2y_0^2}{b'^2x_0\GP{x_0+\cfrac{a}{2}}}}=\cfrac{(a'^2-b'^2)x_0y_0+\cfrac{a}{2}\,a'^2y_0}{b'^2x_0^2+a'^2y_0^2+\cfrac{a}{2}\,b'^2x_0}
\]
Como $a'^2-b'^2=c'^2$ y como $x_0$ y $y_0$, al ser coordenadas de un punto de la elipse, se cumple que $b'^2x_0^2+a'^2y_0^2=a'^2b'^2$, entonces la expresi'on anterior se puede escribir como
\[
\tg{\frac{\al}{2}}=\cfrac{y_0\GP{c'^2x_0+\cfrac{a}{2}\,a'^2}}{b'^2\GP{a'^2+\cfrac{a}{2}\,x_0}}
=\cfrac{y_0\GP{\cfrac{a^2}{4}\,x_0+\cfrac{a}{2}\,a'^2}}{b'^2\GP{a'^2+\frac{a}{2}\,x_0}}=
\frac{a/2}{b'^2}y_0
\]
Como $y_0=b'\sen{E}$ y $b'=\sqrt{a'^2-c'^2}=\sqrt{\GP{\dfrac{b+c}{2}}^2-\GP{\dfrac{a}{2}}^2}$ se tiene finalmente que
\begin{equation}\label{E115}
\tg{\frac{\al}{2}}=\cfrac{a/2}{\sqrt{\GP{\cfrac{b+c}{2}}^2-\GP{\cfrac{a}{2}}^2}}\,\sen{E}
=\cfrac{\cfrac{a}{b+c}}{\sqrt{1-\GP{\cfrac{a}{b+c}}^2}}\,\sen{E}
\end{equation}
De esta ecuaci'on se deduce que, cuando $0\leq E\leq 180^\circ$, 
\[
0\leq\tg{\frac{\al}{2}}\leq\cfrac{\cfrac{a}{b+c}}{\sqrt{1-\GP{\cfrac{a}{b+c}}^2}}
\]
o, en t'erminos de $\sen{\dfrac{\al}{2}}$,
\[
0\leq\sen{\frac{\al}{2}}\leq{\frac{a}{b+c}}
\]
La 'ultima desigualdad da una condici'on para la existencia de una soluci'on en este caso. Sea $\phi$ el 'angulo que forma el semieje menor con uno de los radio vectores en un punto de la elipse tal que $E=90^\circ$ y cuyo seno es $\dfrac{a}{b+c}$. Esto permite escribir la ecuaci'on \eqref{E115} de esta manera:
\begin{equation}\label{E116}
\tg{\frac{\al}{2}}=\tg{\phi}\sen{E}
\end{equation}
Sea $Y_0$ la ordenada del punto Q de la circunferencia auxiliar de la elipse asociado al punto $P$. Dicha ordenada tiene como valor $a'\sen{E}$. Luego,
\begin{equation}\label{E117}
Y_0=a'\cfrac{\tg{\cfrac{\al}{2}}}{\tg{\phi}}\quad\therefore\quad\frac{Y_0}{a'}=\cfrac{\tg{\cfrac{\al}{2}}}{\tg{\phi}}=
\cfrac{\GP{\cfrac{b+c}{2}}\tg{\cfrac{\al}{2}}}{\GP{\cfrac{b+c}{2}}\tg{\phi}}
\end{equation}
Esta proporci'on es la clave para construir el tri'angulo, pues permite buscar el punto en la elipse a traves de su punto correspondiente en su circunferencia auxiliar. Con esto, el proceso para construir el tri'angulo se divide en tres partes a saber:

Primero (Determinaci'on de $Y_0$): Se construye un tri'angulo rect'angulo de hipotenusa $\seg{RS}=a'=\dfrac{b+c}{2}$ y cateto $\seg{ST}=c'=\dfrac{a}{2}$. El 'angulo opuesto a $ST$ ser'a el 'angulo $\phi$ y en su respectivo v'ertice opuesto se traza una circunferencia de radio $\seg{RS}$. Despu'es se traza el 'angulo $A$ con v'ertice en el centro de la circunferencia trazada y uno de sus lados sobre el cateto $RT$ y se biseca. Se traza una tangente a la circunferencia en el punto $D$ donde esta se corta con la prolongaci'on del cateto $RT$ por $T$. Esta tangente intersectar'a a los lados de los 'angulos $\phi$ y $A/2$ en los puntos $E$ y $F$ respectivamente. Los segmentos $\seg{DE}$ y $\seg{DF}$ son a su vez $\dfrac{b+c}{2}\,\tg{\phi}$ y $\dfrac{b+c}{2}\,\tg{\dfrac{\al}{2}}$ (ver figura \ref{F13}-a). Junto con $a'$, se construye el segmento $Y_0$ mediante la cuarta proporcional (figura \ref{F13}-b).

Segundo (Ubicaci'on del punto $Q$ sobre la circunferencia auxiliar y de la abscisa del punto $P$): Se traza un segmento $\seg{BC}=a$. Con centro en la mitad de este segmento $D$ y radio $a'$, se traza una circunferencia. Se traza una perpendicular $p$ a $BC$ por $D$, y a partir de este 'ultimo punto, se busca un punto $E$ sobre $p$ tal que $\seg{DE}=Y_0$. Luego, se traza una paralela a $BC$ por $E$, la cual cortar'a a la circunferencia en un punto $Q$. La longitud del segmento $\seg{EQ}$ corresponde a la abscisa del punto de la circunferencia auxiliar a la elipse en cuesti'on, el cual a su vez es la abscisa del punto para el cual $\al=A$ (figura \ref{F13}-c). 

Tercero (Determinaci'on completa del punto $P$): En la construcci'on anterior, se prolonga el segmento $BC$ por ambos extremos hasta que corte a la circunferencia en los puntos $B'$ y $C'$. Se proyecta el punto $Q$ sobre el segmento $B'C'$, su proyecci'on ser'a el punto $G$. Se ubica un punto $F$ sobre el segmento $DQ$ de tal manera que $\seg{DF}=a'-b'$. Se proyecta el punto $F$ sobre el segmento $B'C'$, dando origen al punto $H$. Se traza una paralela a $DQ$ que pase por $H$, la cual intersectar'a al segmento $GQ$ en $P$. Se unen los puntos $B$ y $C$ con $P$. El tri'angulo $BCP$ es el tri'angulo buscado (Ver figura \ref{F13}-c)\footnote{El m'etodo empleado para trazar el punto de la elipse buscado se basa en tomar una regla y ubicar tres puntos $A$, $B$, $C$, en ella tal que $\seg{AC}=a$ y $\seg{BC}=b$. Luego, se desliza la regla de tal manera que $A$ est'e sobre el eje menor y $B$ est'e sobre el eje mayor, lo cual har'a que $C$ describa la elipse. En \cite{B5} y en \cite{B6}, puede verse una justificaci'on de esta construcci'on.}.
\section{C'alculo del 'area de un tri'angulo oblicu'angulo}\label{S4}
Una vez determinados los par'ametros desconocidos en un tri'angulo oblicu'angulo, es posible determinar otros par'ametros tales como el 'area, per'imetro, bisectrices, etc. En esta secci'on se estudiar'a el c'alculo del 'area.

Se comienza con el 'area de un tri'angulo oblicu'angulo si los datos dados corresponden al cuarto caso, pues a partir de este resultado, se puede encontrar expresiones 'utiles para los casos restantes. Sup'ongase conocidos los lados $a$, $b$ y $C$. Independientemente de si $C$ es agudo u obtuso, el 'area de dicho tri'angulo, de acuerdo a las figuras \ref{F2}-a y \ref{F2}-b es
\[
S=\frac{1}{2}a\cdot \seg{AD}
\]
De estas se deduce que $\seg{AD}=b\sen{C}$, luego
\[
S=\frac{1}{2}ab\sen{C}
\]
Si se trazan alturas respectivas con los lados restantes, un razonamiento similar conduce a lo siguiente:
\begin{equation}\label{E118}
S=\frac{1}{2}ab\sen{C}=\frac{1}{2}bc\sen{A}=\frac{1}{2}ca\sen{B}
\end{equation}
En conclusi'on, {si se conocen dos lados de un tri'angulo y el 'angulo comprendido entre ellos, entonces su 'area es el semiproducto de los lados y el seno del 'angulo ya mencionados}.

Si se toma una ecuaci'on del conjunto de relaciones anteriores, por ejemplo $S=\frac{1}{2}bc\sen{A}$, y si se expresa $b$ y $c$ en funci'on de $a$ mediante el teorema del seno, se tiene que
\[
S=\frac{1}{2}bc\sen{A}=\frac{1}{2}\,\cdot\,\frac{a\sen{B}}{\sen{A}}\,\cdot\,\frac{a\sen{C}}{\sen{A}}\,\cdot\,\sen{A}
=\frac{a^2\sen{B}\sen{C}}{2\sen{A}}
\]
Un proceso similar con las restantes relaciones de 'area, se obtiene el siguiente conjunto de ecuaciones:
\begin{equation}\label{E119}
S=\frac{a^2\sen{B}\sen{C}}{2\sen{A}}=\frac{b^2\sen{A}\sen{C}}{2\sen{B}}=\frac{c^2\sen{A}\sen{B}}{2\sen{C}}
\end{equation}
Est'as ecuaciones permiten calcular el 'area de un tri'angulo si de 'el se conoce un lado y dos 'angulos, ya sea que tengan un lado en com'un o uno de ellos sea opuesto al lado dado. Cuando se conoce del tri'angulo dos lados y un 'angulo opuesto a alguno de ellos, se puede usar ya sea una de las ecuaciones \eqref{E118} o \eqref{E119}. En el caso en el cual se conocen los tres lados, se combina una de las ecuaciones \eqref{E118} con las f'ormulas del 'angulo medio, es decir:
\begin{equation}\label{E120}
S=\frac{1}{2}bc\sen{A}=bc\sen{\frac{A}{2}}\cos{\frac{A}{2}}=bc\sqrt{\frac{p(p-a)}{bc}}\,\sqrt{\frac{(p-b)(p-c)}{bc}}=\sqrt{p(p-a)(p-b)(p-c)}
\end{equation}
Esta ecuaci'on es la muy conocida f'ormula de Her'on. Con ella, se puede calcular el 'area de un tri'angulo conocidos su per'imetro y dos 'angulos internos, una vez calculados sus lados de acuerdo a lo expuesto en la secci'on \ref{S3-6}

El 'area del tri'angulo conocido un lado, su 'angulo opuesto y la suma de los lados restantes es f'acil de calcular. De la figura \ref{F13}-c, y teniendo en cuenta que $a'=\frac{b+c}{2}$, $\sen{\phi}=\frac{a}{b+c}$ y la ecuaci'on \eqref{E117}, se tiene que dicha 'area es
\begin{equation}\label{E121}
S=\frac{1}{2}\seg{BC}\cdot Y_0=\frac{1}{2}a\cdot a'\frac{\tg{\frac{A}{2}}}{\tg{\phi}}
=\frac{1}{4}a(b+c)\tg{\frac{A}{2}}\cotg{\GC{\arcsen{\GP{\frac{a}{b+c}}}}}
\end{equation}
O si expresamos dicha 'area en t'erminos de los par'ametros de la elipse asociada al tri'angulo, entonces
\begin{equation}\label{E122}
S=\frac{1}{4}(2c')(2a')\tg{\frac{A}{2}}\GP{\frac{b'}{c'}}=a'b'\tg{\frac{A}{2}}
\end{equation}
Para esta ecuaci'on, se tuvo en mente que
\[
\tg{\phi}=\frac{\sen{\phi}}{\sqrt{1-\sen^2{\phi}}}=\frac{a}{\sqrt{(b+c)^2-a^2}}=\frac{2c'}{\sqrt{4a'^2-4c'^2}}=\frac{c'}{b'}
\]
\section{Radios de las circunferencias inscrita y circunscrita a un tri'angulo oblicu'angulo}\label{S5}
De las discusiones del teorema del seno y de las f'ormulas del 'angulo medio, surgi'o el radio de la circunferencia circunscrita al tri'angulo en relaci'on con la constante de proporcionalidad del teorema del seno; y el radio de la circunferencia inscrita a un tri'angulo como un medio para simplificar las f'ormulas de las tangentes de los 'angulos medios al igual que en su deducci'on geom'etrica. Adem'as de estas circunferencias, existen otras asociadas a un tri'angulo oblicu'angulo denomidadas excritas, cuya definici'on se da a continuaci'on.

\begin{figure}[htb!]
\centering
\fbox{\includegraphics[scale=0.8]{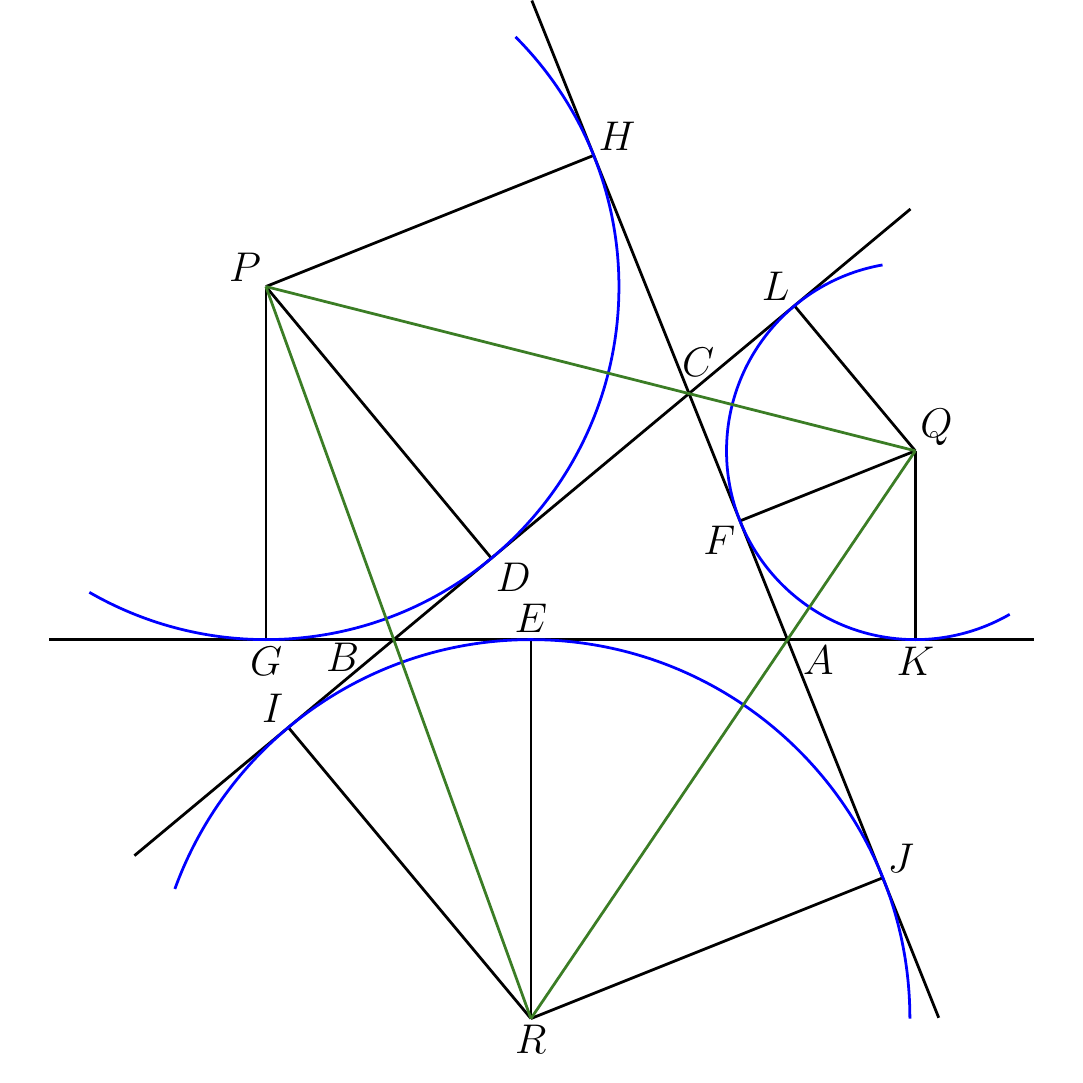}}
\caption{Circunferencias excritas de un tri'angulo oblicu'angulo $ABC$.}\label{F14}
\end{figure}

Se dice que una circunferencia es excrita a un tri'angulo cuando esta es tangente a uno de sus lados y a las prolongaciones de los otros lados tomadas por los extremos del lado en cuesti'on. Consecuencia de esta definici'on es que un tri'angulo cualquiera tiene tres circunferencias excritas (ver figura \ref{F14}).

En la figura mencionada, se unen los centros de las tres circunferencias para formar el tri'angulo $PQR$. Para deducir los radios de estas circunferencias, n'otese que las l'ineas $PQ$, $QR$ y $ST$ son las bisectrices de los tri'angulos externos del tri'angulo $ABC$, pues, como se sabe de la geometr'ia plana, las bisectrices de dos rectas coplanares que se cortan son los lugares geom'etricos de los puntos de ese plano que equidistan de ellas.

Sea $PD=r_a$, $QF=r_b$ y $RE=r_c$ los radios a calcular. Del tri'angulo $BCP$ se conoce el lado $BC=a$ y los 'angulos $\angle PBC=(A+C)/2=90^\circ-B/2$ y $\angle PCB=(A+B)/2=90^\circ-C/2$, par'ametros correspondientes al primer caso de resoluci'on de tri'angulos. Luego, los dem'as par'ametros son iguales a
\begin{align}
\angle{BPC}&=180^\circ-\GP{90^\circ-\frac{B}{2}}-\GP{90^\circ-\frac{C}{2}}=90^\circ-\frac{A}{2}\label{E123}\\
BP&=BC\frac{\sen{\angle{PCB}}}{\sen{\angle{BPC}}}=a\cfrac{\cos{\cfrac{C}{2}}}{\cos{\cfrac{A}{2}}}\label{E124}\\
CP&=BC\frac{\sen{\angle{PBC}}}{\sen{\angle{BPC}}}=a\cfrac{\cos{\cfrac{B}{2}}}{\cos{\cfrac{A}{2}}}\label{E125}
\end{align}
El radio de la circunferencia excrita de centro $P$ y radio $r_a$ es igual a 
\begin{equation}\label{E126}
r_a=PD=PB\sen{\angle{PBC}}=a\cfrac{\cos{\cfrac{B}{2}}\cos{\cfrac{C}{2}}}{\cos{\cfrac{A}{2}}}
\end{equation}
De las f'ormulas del 'angulo medio, el numerador de la ecuaci'on anterior se puede escribir de la siguiente manera:
\begin{equation}\label{E127}
\cos{\frac{B}{2}}\cos{\frac{C}{2}}=\sqrt{\frac{p^2(p-b)(p-c)}{a^2bc}}=\frac{p}{a}\sen{\frac{A}{2}}
\end{equation}
Este relaci'on permite escribir el radio en cuesti'on de dos maneras:
\begin{equation}\label{E128}
r_a=p\tg{\frac{A}{2}}=\frac{pr}{p-a}
\end{equation}
La 'ultima relaci'on surge al usar las formulas de las tangentes de los 'angulos medios.
Los radios restantes se obtienen con los par'ametros resueltos de los tri'angulos $AQC$ y $BRC$.

Para el tri'angulo $AQC$: se conocen el lado $AC=b$ y los 'angulos $\angle{QAC}=90^\circ-\frac{A}{2}$ y $\angle{QCA}=90^\circ-\frac{C}{2}$, sus dem'as par'ametros son:
\begin{align}
\angle{AQC}&=180^\circ-\angle{QAC}-\angle{QCA}\notag\\
&=180^\circ-\GP{90^\circ-\frac{A}{2}}-\GP{90^\circ-\frac{C}{2}}=90^\circ-\frac{B}{2}\label{E129}\\
AQ&=AC\frac{\sen{\angle{QCA}}}{\sen{\angle{AQC}}}=b\cfrac{\cos{\cfrac{C}{2}}}{\cos{\cfrac{B}{2}}}\label{E130}\\
CQ&=AC\frac{\sen{\angle{QAC}}}{\sen{\angle{AQC}}}=b\cfrac{\cos{\cfrac{A}{2}}}{\cos{\cfrac{B}{2}}}\label{E131}
\end{align}
El radio $r_b$ es igual a 
\begin{equation}\label{E132}
r_b=QF=CQ\sen{\angle{QCA}}=b\cfrac{\cos{\cfrac{A}{2}}\cos{\cfrac{C}{2}}}{\cos{\cfrac{B}{2}}}
\end{equation}
y combinando esta relaci'on con las f'ormulas del 'angulo medio, se tiene que 
\begin{equation}\label{E133}
r_b=p\tg{\frac{B}{2}}=\frac{pr}{p-b}
\end{equation}
Para el tri'angulo $BRC$: se conocen el lado $AB=c$ y los 'angulos $\angle{RAB}=90^\circ-\frac{A}{2}$ y $\angle{RBA}=90^\circ-\frac{B}{2}$, sus dem'as par'ametros son:
\begin{align}
\angle{ARB}&=180^\circ-\angle{RAB}-\angle{RBA}\notag\\
&=180^\circ-\GP{90^\circ-\frac{A}{2}}-\GP{90^\circ-\frac{B}{2}}=90^\circ-\frac{C}{2}\label{E134}\\
AR&=AB\frac{\sen{\angle{RBA}}}{\sen{\angle{ARB}}}=c\cfrac{\cos{\cfrac{B}{2}}}{\cos{\cfrac{C}{2}}}\label{E135}\\
BR&=AB\frac{\sen{\angle{RAB}}}{\sen{\angle{ARB}}}=c\cfrac{\cos{\cfrac{A}{2}}}{\cos{\cfrac{C}{2}}}\label{E136}
\end{align}
El radio $r_c$ es igual a 
\begin{equation}\label{E137}
r_c=RE=BR\sen{\angle{RBA}}=c\cfrac{\cos{\cfrac{A}{2}}\cos{\cfrac{B}{2}}}{\cos{\cfrac{C}{2}}}
\end{equation}
y combinando esta relaci'on con las f'ormulas del 'angulo medio, se tiene que 
\begin{equation}\label{E138}
r_c=p\tg{\frac{C}{2}}=\frac{pr}{p-c}
\end{equation}
Sumando los inversos de los radios de las circunferencias excritas, se tiene que
\begin{equation}\label{E139}
\frac{1}{r_a}+\frac{1}{r_b}+\frac{1}{r_c}=\frac{p-a}{pr}+\frac{p-b}{pr}+\frac{p-c}{pr}=\frac{1}{r}
\end{equation}
Esto muestra que el radio de la circunferencia circunscrita a un tri'angulo oblicu'angulo es media arm'onica de los radios de las circunferencias excritas a dicho tri'angulo.

En la secci'on \ref{S2}-\ref{S2-1} se mostro que el di'ametro de la circunferencia circunscrita a un tri'angulo es la constante de proporcionalidad del teorema del seno:
\begin{equation}\label{E140}
\frac{a}{\sen{A}}=\frac{b}{\sen{B}}=\frac{c}{\sen{C}}=2R
\end{equation}
Combinando esta relaci'on con las ecuaciones \eqref{E126}, \eqref{E132} y \eqref{E137}, se obtienen los radios de las circunferencias excritas en funci'on del radio de la circunferencia circunscrita y los 'angulos internos del tri'angulo:
\begin{align}
r_a&=a\cfrac{\cos{\cfrac{B}{2}}\cos{\cfrac{C}{2}}}{\cos{\cfrac{A}{2}}}=
2R\sen{A}\cfrac{\cos{\cfrac{B}{2}}\cos{\cfrac{C}{2}}}{\cos{\cfrac{A}{2}}}
=4R\sen{\cfrac{A}{2}}\cos{\cfrac{B}{2}}\cos{\cfrac{C}{2}}\label{E141}\\
r_b&=b\cfrac{\cos{\cfrac{A}{2}}\cos{\cfrac{C}{2}}}{\cos{\cfrac{B}{2}}}=
2R\sen{B}\cfrac{\cos{\cfrac{A}{2}}\cos{\cfrac{C}{2}}}{\cos{\cfrac{B}{2}}}=
4R\cos{\cfrac{A}{2}}\sen{\cfrac{B}{2}}\cos{\cfrac{C}{2}}\label{E142}\\
r_c&=c\cfrac{\cos{\cfrac{A}{2}}\cos{\cfrac{B}{2}}}{\cos{\cfrac{C}{2}}}=
2R\sen{C}\cfrac{\cos{\cfrac{A}{2}}\cos{\cfrac{B}{2}}}{\cos{\cfrac{C}{2}}}=
4R\cos{\cfrac{A}{2}}\cos{\cfrac{B}{2}}\sen{\cfrac{C}{2}}\label{E143}
\end{align}
Reemplazando las ecuaciones anteriores en la ecuaci'on \eqref{E139}:
\begin{align*}
\frac{1}{r}&=\frac{1}{r_a}+\frac{1}{r_b}+\frac{1}{r_c}\\
&=\cfrac{1}{4R\sen{\cfrac{A}{2}}\cos{\cfrac{B}{2}}\cos{\cfrac{C}{2}}}
+\cfrac{1}{4R\cos{\cfrac{A}{2}}\sen{\cfrac{B}{2}}\cos{\cfrac{C}{2}}}
+\cfrac{1}{4R\cos{\cfrac{A}{2}}\cos{\cfrac{B}{2}}\sen{\cfrac{C}{2}}}\\
&=\cfrac{\cos{\cfrac{A}{2}}\sen{\cfrac{B}{2}}\sen{\cfrac{C}{2}}+\sen{\cfrac{A}{2}}\cos{\cfrac{B}{2}\sen{\cfrac{C}{2}}+\sen{\cfrac{A}{2}}\sen{\cfrac{B}{2}}}\cos{\cfrac{C}{2}}}{4R\sen{\cfrac{A}{2}}\cos{\cfrac{A}{2}}\sen{\cfrac{B}{2}}\cos{\cfrac{B}{2}}\sen{\cfrac{C}{2}}\cos{\cfrac{C}{2}}}\\
&=\cfrac{\sen{\cfrac{A+B}{2}}\sen{\cfrac{C}{2}}+\sen{\cfrac{A}{2}}\sen{\cfrac{B}{2}}\cos{\cfrac{C}{2}}}{4R\sen{\cfrac{A}{2}}\cos{\cfrac{A}{2}}\sen{\cfrac{B}{2}}\cos{\cfrac{B}{2}}\sen{\cfrac{C}{2}}\cos{\cfrac{C}{2}}}
\end{align*}
Como $A+B+C=180^\circ$, $\sen{\dfrac{A+B}{2}}=\cos{\dfrac{C}{2}}$ y $\sen{\dfrac{C}{2}}=\cos{\dfrac{A+B}{2}}$. Al colocar estas relaciones en la ecuaci'on anterior, se tiene que
\begin{align*}
\frac{1}{r}&=\cfrac{\sen{\cfrac{A+B}{2}}\sen{\cfrac{C}{2}}+\sen{\cfrac{A}{2}}\sen{\cfrac{B}{2}}\cos{\cfrac{C}{2}}}{4R\sen{\cfrac{A}{2}}\cos{\cfrac{A}{2}}\sen{\cfrac{B}{2}}\cos{\cfrac{B}{2}}\sen{\cfrac{C}{2}}\cos{\cfrac{C}{2}}}\\
&=\cfrac{\cos{\cfrac{C}{2}}\NC{\cos{\cfrac{A+B}{2}}+\sen{\cfrac{A}{2}}\sen{\cfrac{B}{2}}}}{4R\sen{\cfrac{A}{2}}\cos{\cfrac{A}{2}}\sen{\cfrac{B}{2}}\cos{\cfrac{B}{2}}\sen{\cfrac{C}{2}}\cos{\cfrac{C}{2}}}\\
&=\cfrac{1}{4R\sen{\cfrac{A}{2}}\sen{\cfrac{B}{2}}\sen{\cfrac{C}{2}}}
\end{align*}
La relaci'on obtenida, al reescribirla de la siguiente manera:
\begin{equation}\label{E144}
r=4R\sen{\frac{A}{2}}\sen{\frac{B}{2}}\sen{\frac{C}{2}}
\end{equation}
ilustra una relaci'on entre los radios de las circunferencias inscrita y circunscrita en un tri'angulo oblicu'angulo. Una consecuencia de las ecuaciones \eqref{E141}--\eqref{E144} es que el radio de la circunferencia inscrita al igual que los radios de las circunferencias excritas son menores que el doble del di'ametro de la circunferencia circunscrita.

El tri'angulo $PQR$, el cual nace al unir los centros de las circunferencias excritas, se denomina el tri'angulo excrito asociado al tri'angulo oblicu'angulo $ABC$. Sus 'angulos internos se determinaron como parte del procedimiento en la deducci'on de los radios de las circunferencias en cuesti'on. Sus valores son:
\begin{align}
P&\equiv\angle{BPC}=90^\circ-\frac{A}{2}\label{E145}\\
Q&\equiv\angle{CQA}=90^\circ-\frac{B}{2}\label{E146}\\
R&\equiv\angle{ARB}=90^\circ-\frac{C}{2}\label{E147}
\end{align}
Solo resta determinar sus lados. El lado $PQ$ es igual a:
\[
PQ=PC+CQ=a\cfrac{\cos{\cfrac{B}{2}}}{\cos{\cfrac{A}{2}}}+b\cfrac{\cos{\cfrac{A}{2}}}{\cos{\cfrac{B}{2}}}
\]
Expresando $a$ y $b$ en funci'on de $c$ mediante el teorema del seno:
\begin{align}
PQ&=\frac{c}{\sen{C}}\GP{\frac{\sen{A}\cos{\frac{B}{2}}}{\cos{\frac{A}{2}}}+\frac{\sen{B}\cos{\frac{A}{2}}}{\cos{\frac{B}{2}}}}\notag\\
&=\frac{2c}{\sen{C}}\sen{\frac{A+B}{2}}=\frac{2c}{\sen{C}}\cos{\frac{C}{2}}=\cfrac{c}{\sen{\cfrac{C}{2}}}\label{E148}
\end{align}
Mediante un procedimiento similar, los lados $QR$ y $RP$ son iguales a
\begin{align}
QR&=QA+AR=b\cfrac{\cos{\cfrac{C}{2}}}{\cos{\cfrac{B}{2}}}+c\cfrac{\cos{\cfrac{B}{2}}}{\cos{\cfrac{C}{2}}}=\cfrac{a}{\sen{\cfrac{A}{2}}}\label{E149}\\
RP&=RB+BP=c\cfrac{\cos{\cfrac{A}{2}}}{\cos{\cfrac{C}{2}}}+a\cfrac{\cos{\cfrac{C}{2}}}{\cos{\cfrac{A}{2}}}=\cfrac{b}{\sen{\cfrac{B}{2}}}\label{E150}
\end{align}
Otra manera de obtener los lados del tri'angulo excrito $PQR$ es mediante la semejanza que este tiene con los tri'angulos $BCP$, $ACQ$ y $ABR$, los cuales a su vez, son semejantes entre s'i. Solo basta resolver uno de estos tri'angulos para resolver los dem'as.

Para finalizar esta secci'on, se discutir'a el 'area de un tri'angulo oblicu'angulo en relaci'on con los radios de sus circunferencias asociadas. Al despejar los senos de los 'angulos internos en las ecuaciones \eqref{E140}:
\begin{align}
\sen{A}&=\frac{a}{2R}\label{E151}\\
\sen{B}&=\frac{b}{2R}\label{E152}\\
\sen{C}&=\frac{c}{2R}\label{E153}
\end{align}
y reemplaz'andolos en las ecuaciones \eqref{E118}, se tiene que
\begin{align*}
S&=\frac{1}{2}ab\sen{C}=\frac{ab}{2}\GP{\frac{c}{2R}}=\frac{abc}{4R}\\
S&=\frac{1}{2}bc\sen{A}=\frac{bc}{2}\GP{\frac{a}{2R}}=\frac{abc}{4R}\\
S&=\frac{1}{2}ca\sen{B}=\frac{ca}{2}\GP{\frac{b}{2R}}=\frac{abc}{4R}
\end{align*}
El 'area de un tri'angulo en t'erminos del radio de la circunferencia circunscrita es
\begin{equation}\label{E154}
S=\frac{abc}{4R}
\end{equation}
Si se escribe la f'ormula de Her'on de la siguiente manera:
\[
S=\sqrt{p(p-a)(p-b)(p-c)}=p\sqrt{\frac{(p-a)(p-b)(p-c)}{p}}
\]
el t'ermino que acompa'na a $p$ es el radio de la circunferencia inscrita, lo que permite expresar el 'area del tri'angulo en funci'on de dicho radio:
\begin{equation}\label{E155}
S=pr
\end{equation}
Si se multiplican los cosenos de los 'angulos medios:
\begin{align*}
\cos{\frac{A}{2}}\cos{\frac{B}{2}}\cos{\frac{C}{2}}&=
\sqrt{\frac{p(p-a)}{bc}}\,\sqrt{\frac{p(p-b)}{ac}}\,\sqrt{\frac{p(p-c)}{ab}}\\&=
\frac{p}{abc}\sqrt{p(p-a)(p-b)(p-c)}
\end{align*}
Utilizando las ecuaciones \eqref{E154} y \eqref{E155} en la relaci'on anterior:
\begin{equation}\label{E156}
\cos{\frac{A}{2}}\cos{\frac{B}{2}}\cos{\frac{C}{2}}=\frac{p}{S/4R}\,S=\frac{p}{4R}
\end{equation}
Por otro lado, al multiplicar las ecuaciones \eqref{E126}, \eqref{E132} y \eqref{E137} miembro a miembro:
\[
r_ar_br_c=abc\cos{\frac{A}{2}}\cos{\frac{B}{2}}\cos{\frac{C}{2}}
\]
Y al combinarla con la ecuaci'on \eqref{E156}:
\[
r_ar_br_c=p\,\frac{abc}{4R}
\]
De las ecuaciones \eqref{E154} y \eqref{E155}, la ecuaci'on anterior se convierte en $r_ar_br_c=S^2/r$, o despejando el cuadrado del 'area:
\begin{equation}\label{E157}
S^2=rr_ar_br_c
\end{equation}
\section{Medianas, bisectrices y alturas de un tri'angulo oblicu'angulo}\label{S6}
En esta secci'on, se estudiar'a las relaciones entre los lados y 'angulos internos de un tri'angulo oblicu'angulo con sus medianas, bisectrices y alturas.
\subsection{Medianas}\label{S6-1} 
\begin{figure}[htb!]
\centering
\fbox{\includegraphics[scale=0.8]{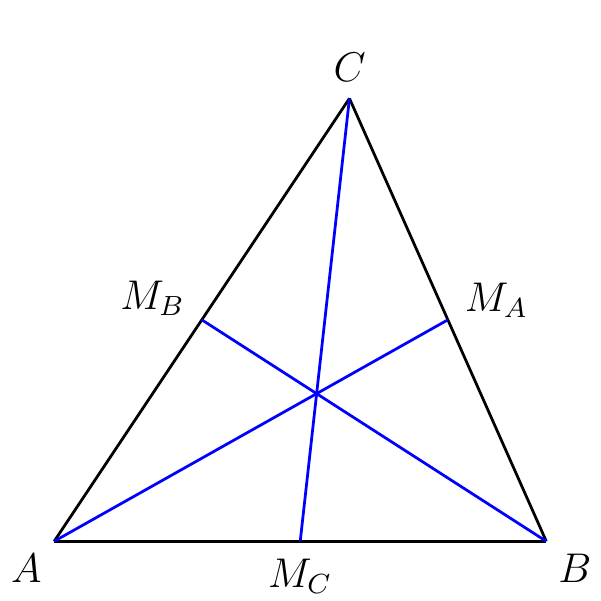}}
\caption{Medianas de un tri'angulo oblicu'angulo}\label{F15}
\end{figure}

La mediana de un tri'angulo es un segmento que une un v'ertice cualquiera con el punto medio de su lado opuesto (ver figura \ref{F15}). Entre las propiedades conocidas de las medianas est'an que se intersectan en un punto denominado baricentro, y que dicho punto dista de un v'ertice cualquiera $2/3$ de su mediana correspondiente.

Sea $M_A$, $M_B$ y $M_C$ los puntos medios de los lados $BC$, $CA$ y $AB$ respectivamente y sea $AM_A$, $BM_B$ y $CM_C$ las medianas del tri'angulo oblicu'angulo $ABC$, cuyas longitudes correspondiente son $m_a$, $m_b$ y $m_c$. (figura \ref{F15}). Para calcular la longitud de una de ellas, por ejemplo $m_a$, se usa uno de los tri'angulos $ABM_A$ o $ACM_A$. Ambos tri'angulos tri'angulos pertenecen al cuarto caso, pues se conocen dos lados y un 'angulo comprendido entre ellos.

Aplicando el teorema del coseno al tri'angulo $ABM_A$:
\[
m_a^2=\overline{AM_A}^2=
{\overline{AB}^2+\overline{BM_B}^2-2\overline{AB}\cdot\overline{BM_B}\cos{B}}=
{c^2+\frac{a^2}{4}-ac\cos{B}}
\]
y haciendo lo mismo para el tri'angulo $ACM_A$:
\[
m_a^2=\overline{AM_A}^2=
{\overline{AC}^2+\overline{CM_C}^2-2\overline{AC}\cdot\overline{CM_C}\cos{C}}
={b^2+\frac{a^2}{4}-ab\cos{C}}
\]
Por lo tanto
\begin{equation}\label{E158}
m_a^2={c^2+\frac{a^2}{4}-ac\cos{B}}={b^2+\frac{a^2}{4}-ab\cos{C}}
\end{equation}
Un an'alisis similar permite establecer relaciones similares con las medianas restantes:
\begin{align}
m_b^2&={c^2+\frac{b^2}{4}-bc\cos{A}}={a^2+\frac{b^2}{4}-ab\cos{C}}\label{E159}\\
m_c^2&={b^2+\frac{c^2}{4}-bc\cos{A}}={a^2+\frac{c^2}{4}-ac\cos{B}}\label{E160}
\end{align}
Si se expresa los cosenos de los 'angulos internos en funci'on de sus lados mediante el teorema del coseno, las ecuaciones \eqref{E158}--\eqref{E160} se convierten en 
\begin{align}
m_a^2=&\frac{b^2+c^2}{2}-\frac{a^2}{4}\label{E161}\\
m_b^2=&\frac{a^2+c^2}{2}-\frac{b^2}{4}\label{E162}\\
m_c^2=&\frac{a^2+b^2}{2}-\frac{c^2}{4}\label{E163}
\end{align}
De esta manera, se pueden expresar los cuadrados de las medianas en funci'on de los lados. 

Un caso de resoluci'on de tri'angulos es determinar los lados y 'angulos de un tri'angulo oblicu'angulo si se conocen las longitudes de sus medianas. En ese caso, las ecuaciones anteriores son 'utiles para ello. Estas forman un sistema lineal de tres ecuaciones con tres inc'ognitas, las cuales son los cuadrados de los lados. Al resolver este sistema, se tiene que 
\begin{align}
\frac{9}{4}a^2&=2(m_b^2+m_c^2)-m_a^2\label{E164}\\
\frac{9}{4}b^2&=2(m_a^2+m_c^2)-m_b^2\label{E165}\\
\frac{9}{4}c^2&=2(m_a^2+m_b^2)-m_c^2\label{E166}
\end{align}
De estas ecuaciones, se observa que no todos los valores de las medianas dadas pertenecen a un tri'angulo. Para que esto suceda, los lados derechos de dichas ecuaciones deben ser mayores que cero, es decir:
\begin{align}
m_b^2+m_c^2&>\frac{m_a^2}{2}\label{E167}\\
m_a^2+m_c^2&>\frac{m_b^2}{2}\label{E168}\\
m_a^2+m_b^2&>\frac{m_c^2}{2}\label{E169}
\end{align}
Si las medianas dadas cumplen las desigualdades anteriores, los lados del tri'angulo obtenidos con las ecuaciones \eqref{E164}--\eqref{E166} deben satisfacer las condiciones deducidas en el quinto caso de resoluci'on de tri'angulos, ya sea las que est'an dadas por las desigualdades \eqref{E75}--\eqref{E77}, o que el semiper'imetro sea mayor que los tres lados. Para finalizar esta secci'on, se encontrar'a el 'angulo que forma una mediana con su lado del tri'angulo correspondiente, en particular, el 'angulo $\angle AM_CC$. Para ello, se aplica el teorema del seno y del coseno al tri'angulo $AM_CC$:
\begin{align}
\sen{\angle AM_CC}&=\frac{b\sen{A}}{m_c}\label{E170}\\
\cos{\angle AM_CC}&=\frac{c^2/4+m_c^2-b^2}{cm_c}\label{E171}
\end{align}
Dividiendo las ecuaciones anteriores miembro a miembro y combinando este resultado con la ecuaci'on \eqref{E163}:
\begin{equation}\label{E172}
\cotg{\angle AM_CC}=\frac{c^2/4+m_c^2-b^2}{bc\sen{A}}=\frac{a^2-b^2}{2bc\sen{A}}
\end{equation}
De las relaciones de Mollweide, se puede expresar la diferencia de cuadrados de $a$ y $b$, multiplicando miembro a miembro las ecuaciones \eqref{E42}:
\[
\frac{a+b}{c}\,\frac{a-b}{c}=\frac{a^2-b^2}{c^2}=\cfrac{\cos{{\cfrac{A-B}{2}}}}{\sen{\cfrac{C}{2}}}\,\cfrac{\sen{{\cfrac{A-B}{2}}}}{\cos{\cfrac{C}{2}}}
=\frac{\sen{(A-B)}}{\sen{C}}
\]
Despejando la diferencia en cuesti'on y reemplaz'andola en la ecuaci'on \eqref{E175}, se tiene que
\begin{equation}\label{E173}
\cotg{\angle AM_CC}=\frac{c}{2b\sen{A}}\frac{\sen{(A-B)}}{\sen{C}}
\end{equation}
Del teorema del seno, el cociente $c/b$ es igual a $\sen{C}/\sen{B}$, luego
\begin{equation}\label{E174}
\cotg{\angle AM_CC}=\frac{\sen{(A-B)}}{2\sen{A}\sen{B}}=\frac{\cotg{B}-\cotg{A}}{2}
\end{equation}
El 'angulo $\angle AM_CC$ ser'a agudo si $B<A$ y ser'a obtuso si $B>A$. El 'angulo $BM_CC$ es el suplemento del 'angulo $\angle AM_CC$. Los 'angulos que forman las medianas restantes con sus respectivos lados se obtienen de manera similar.
\subsection{Bisectrices}\label{S6-2}
\begin{figure}[htb!]
\centering
\fbox{\includegraphics[scale=0.8]{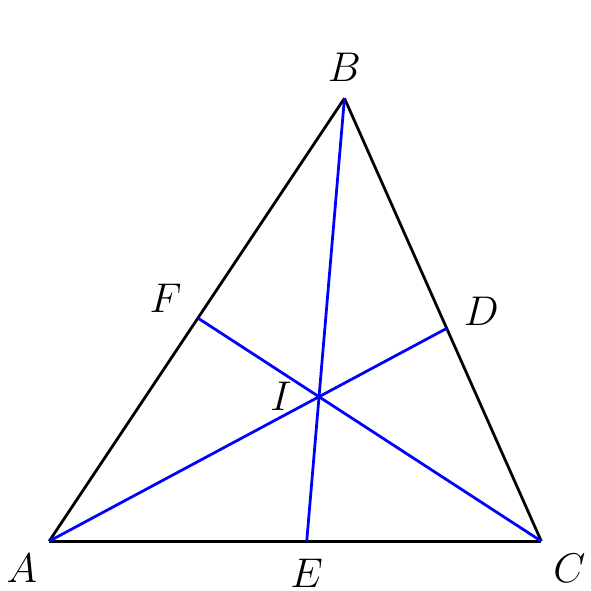}}
\caption{Bisectrices de un tri'angulo oblicu'angulo}\label{F16}
\end{figure}

La bisectriz de un tri'angulo oblicu'angulo es una l'inea recta que biseca un 'angulo interno de un tri'angulo.  Consecuencia de esta definici'on es que un tri'angulo cualquiera tiene tres de ellas. De la geometr'ia elemental, se sabe que se intersectan en un punto $I$ denominado incentro, el cual es el centro de la circunferencia inscrita al tri'angulo (ver figura \ref{F16}).

Una propiedad muy conocida de las bisectrices de un tri'angulo es que una de ellas divide al lado opuesto en segmentos que son proporcionales a los lados restantes. Este resultado puede probarse como sigue: Aplicando el teorema del seno al tri'angulo $ABD$:
\[
\frac{\seg{BD}}{\sen{\angle{BAD}}}=\frac{\seg{AB}}{\sen{\angle{ADB}}} 
\]
Haciendo lo mismo, pero con el tri'angulo $ACD$, se tiene que
\[
\frac{\seg{DC}}{\sen{\angle{CAD}}}=\frac{\seg{AC}}{\sen{\angle{ADC}}}
\]
Como $AD$ es la bisectriz del 'angulo $\angle{BAC}$, $\angle{BAD}=\angle{CAD}$; y como los 'angulos $\angle{ADB}$ y $\angle{ADC}$ son suplementarios, $\sen{\angle{ADB}}=\sen{(180^\circ-\angle{ADC})}=\sen{\angle{ADC}}$, entonces al dividir miembro a miembro las proporciones anteriores y teniendo en cuenta estas observaciones, se deduce que
\begin{equation}\label{E175}
\frac{\,\seg{BD}\,}{\,\seg{DC}\,}=\frac{\,\seg{AB}\,}{\,\seg{AC}\,}
\end{equation}
De esta relaci'on se deducen las siguientes dos\footnote{Si $\dfrac{a}{b}=\dfrac{c}{d}$, entonces
\begin{align*}
\frac{a}{a+b}&=\frac{c}{c+d}\\
\frac{a+b}{b}&=\frac{c+d}{d}
\end{align*}}:
\begin{align}
\frac{\seg{BD}}{\seg{BD}+\seg{DC}}&=\frac{\,\seg{BD}\,}{\,\seg{BC}\,}=\frac{\seg{AB}}{\seg{AB}+\seg{AC}}\quad\therefore\quad \seg{BD}=\frac{ac}{b+c}\label{E176}\\
\frac{\seg{BD}+\seg{DC}}{\seg{DC}}&=\frac{\,\seg{BC}\,}{\,\seg{DC}\,}=\frac{\seg{AB}+\seg{AC}}{\seg{AC}}\quad\therefore\quad \seg{DC}=\frac{ab}{b+c}\label{E177}
\end{align}
Aplicando el teorema del seno al tri'angulo $ABD$:
\begin{equation}\label{E178}
\frac{\seg{AD}}{\sen{B}}=\cfrac{\seg{BD}}{\sen{\cfrac{A}{2}}}\quad\therefore\quad
\seg{AD}=\seg{BD}\cfrac{\sen{B}}{\sen{\cfrac{A}{2}}}
\end{equation}
Como $\sen{B}=\dfrac{b}{a}\sen{A}$, entonces al reemplazar esta relaci'on y la ecuaci'on \eqref{E179} en la ecuaci'on \eqref{E181}, se deduce que la longitud de la bisectriz $AD$ es
\begin{equation}\label{E179}
\seg{AD}=\frac{ac}{b+c}\cdot\frac{b}{a}\cdot\cfrac{\sen{A}}{\sen{\cfrac{A}{2}}}
=\frac{2bc}{b+c}\cos{\cfrac{A}{2}}
\end{equation}
Para expresar esta bisectriz en funci'on de los lados, se reemplaza la ecuaci'on \eqref{E22} correspondiente para el coseno del 'angulo medio, es decir:
\begin{equation}\label{E180}
\seg{AD}=\frac{2bc}{b+c}\sqrt{\frac{p(p-a)}{bc}}=\sqrt{bc\cdot\frac{(a+b+c)(b+c-a)}{(b+c)^{2}}}=
\sqrt{bc}\,\sqrt{1-\GP{\frac{a}{b+c}}^{2}}
\end{equation}
Un razonamiento similar conducir'a a relaciones an'alogas para las dem'as bisectrices. Estas se enuncian a continuaci'on:

Para la bisectriz $BE$:
\begin{align}
\seg{AE}&=\frac{bc}{a+c}\label{E181}\\
\seg{EC}&=\frac{ab}{a+c}\label{E182}\\
\seg{BE}&=\frac{2ac}{a+c}\cos{\frac{B}{2}}=\sqrt{ac}\,\sqrt{1-\GP{\frac{b}{a+c}}^2}\label{E183}
\end{align}
y para la bisectriz $CF$:
\begin{align}
\seg{AF}&=\frac{bc}{a+b}\label{E184}\\
\seg{FB}&=\frac{ac}{a+b}\label{E185}\\
\seg{CF}&=\frac{2ab}{a+b}\cos{\frac{C}{2}}=\sqrt{ab}\,\sqrt{1-\GP{\frac{c}{a+b}}^2}\label{E186}
\end{align}
Las relaciones anteriormente deducidas no son las 'unicas que pueden derivarse para determinar las longitudes de las bisectrices. Para encontrarlas, se determinar'an los 'angulos que forma la bisectriz $AD$ con el lado $BC$ en funci'on de los 'angulos adyacentes a estos lados. Los 'angulos $\angle ADC$ y $\angle ADB$ son exteriores a los tri'angulos $ADC$ y $ADB$, respectivamente, luego
\begin{align*}
\angle ADC&=\frac{A}{2}+B\\
\angle ADB&=\frac{A}{2}+C
\end{align*} 
Como $A=180^\circ-B-C$, estos 'angulos se convierten en 
\begin{align}
\angle ADC&=90^\circ+\frac{B-C}{2}\label{E187}\\
\angle ADB&=90^\circ-\frac{B-C}{2}\label{E188}
\end{align}
Aplicando el teorema del seno tanto al tri'angulo $ADC$ como al tri'angulo $ADB$:
\begin{align*}
\frac{\seg{AD}}{\sen{B}}&=\cfrac{c}{\sen{\GP{90^\circ-\cfrac{B-C}{2}}}}
=\cfrac{c}{\cos{\cfrac{B-C}{2}}}\\
\cfrac{\seg{AD}}{\sen{C}}&=\cfrac{b}{\sen{\GP{90^\circ-\cfrac{B-C}{2}}}}=\cfrac{b}{\cos{\cfrac{B-C}{2}}}
\end{align*}   
De estas ecuaciones se obtiene la bisectriz $AD$. Un razonamiento similar concudir'a a relaciones similares para las dem'as bisectrices, las cuales se presentan a continuaci'on:
\begin{align}
\seg{AD}&=\cfrac{b\sen{C}}{\cos{\cfrac{B-C}{2}}}=\cfrac{c\sen{B}}{\cos{\cfrac{B-C}{2}}}\label{E189}\\
\seg{BE}&=\cfrac{c\sen{A}}{\cos{\cfrac{C-A}{2}}}=\cfrac{a\sen{C}}{\cos{\cfrac{C-A}{2}}}\label{E190}\\
\seg{CF}&=\cfrac{a\cos{B}}{\cos{\cfrac{A-B}{2}}}=\cfrac{b\cos{A}}{\cos{\cfrac{A-B}{2}}}\label{E191}
\end{align}
Para finalizar esta secci'on, se calcular'a la distancia del incentro a los v'ertices de un tri'angulo oblicu'angulo y la raz'on entre esta distancia y la longitud de su bisectriz correspondiente. Para ello, el tri'angulo $AIB$ ser'a de utilidad. De 'el se conoce dos 'angulos y un lado com'un a ellos. De las relaciones del caso I, se deduce que $\seg{AI}$ es igual a:
\begin{equation}\label{E192}
\seg{AI}=\cfrac{c\sen{\cfrac{B}{2}}}{\sen{\GP{180^\circ-\cfrac{A+B}{2}}}}=
\cfrac{c\sen{\cfrac{B}{2}}}{\cos{\cfrac{C}{2}}}
\end{equation}
Si se razona de igual manera usando el tri'angulo $AIC$, esta distancia es igual a:
\begin{equation}\label{E193}
\seg{AI}=\cfrac{b\sen{\cfrac{C}{2}}}{\cos{\cfrac{B}{2}}}
\end{equation}
Si se dividen las ecuaciones \eqref{E195} y \eqref{E196} con las ecuaciones \eqref{E192} de tal manera que tanto los terminos $b$ o $c$ se eliminen, se obtiene lo siguiente:
\begin{equation}\label{E194}
\frac{\,\seg{AI}\,}{\,\seg{AD}\,}=\cfrac{\cos{\cfrac{B-C}{2}}}{2\cos{\cfrac{B}{2}}\cos{\cfrac{C}{2}}}
=\frac{1}{2}\GP{1+\tg{\frac{B}{2}}\tg{\frac{C}{2}}}
\end{equation}
El producto de las tangentes en esta ecuaci'on se calcul'o en la secci'on \ref{S3-6}, luego, al reemplazar la ecuaci'on \eqref{E88} en la ecuaci'on \eqref{E197}, se tendr'a que
\begin{equation}\label{E195}
\frac{\,\seg{AI}\,}{\,\seg{AD}\,}=1-\frac{a}{2p}=\frac{b+c}{a+b+c}
\end{equation}
De forma similar se puede calcular las distancias del incentro a los dem'as v'ertices y sus razones con sus correspondientes longitudes de sus bisectrices. Estas se presentan a continuaci'on:

Para $\seg{BI}$:
\begin{align}
\seg{BI}&=\cfrac{a\sen{\cfrac{C}{2}}}{\cos{\cfrac{A}{2}}}=\cfrac{c\sen{\cfrac{A}{2}}}{\cos{\cfrac{C}{2}}}\label{E196}\\
\frac{\,\seg{BI}\,}{\,\seg{BE}\,}&=\cfrac{\cos{\cfrac{C-A}{2}}}{2\cos{\cfrac{C}{2}}\cos{\cfrac{A}{2}}}
=\frac{1}{2}\GP{1+\tg{\frac{C}{2}}\tg{\frac{A}{2}}}=1-\frac{b}{2p}=\frac{a+c}{a+b+c}\label{E197}
\end{align}
Para $\seg{CI}$:
\begin{align}
\seg{CI}&=\cfrac{a\sen{\cfrac{B}{2}}}{\cos{\cfrac{A}{2}}}=\cfrac{b\sen{\cfrac{A}{2}}}{\cos{\cfrac{B}{2}}}\label{E198}\\
\frac{\,\seg{CI}\,}{\,\seg{CE}\,}&=\cfrac{\cos{\cfrac{A-B}{2}}}{2\cos{\cfrac{A}{2}}\cos{\cfrac{B}{2}}}
=\frac{1}{2}\GP{1+\tg{\frac{A}{2}}\tg{\frac{B}{2}}}=1-\frac{c}{2p}=\frac{a+b}{a+b+c}\label{E199}
\end{align}
Es interesante resaltar que la suma de las tres razones da como resultado:
\begin{equation}\label{E200}
\frac{\,\seg{AI}\,}{\,\seg{AD}\,}+\frac{\,\seg{BI}\,}{\,\seg{BE}\,}+\frac{\,\seg{CI}\,}{\,\seg{CE}\,}=2
\end{equation}
\subsection{Alturas}\label{S6-3}
\begin{figure}[htb!]
\centering
\begin{tabular}{cccc}
\fbox{\includegraphics[scale=0.8]{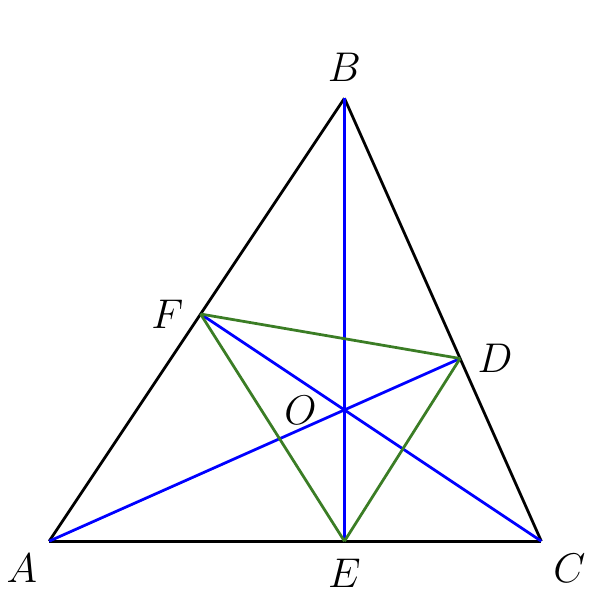}}\\
a)\\
\fbox{\includegraphics[scale=0.8]{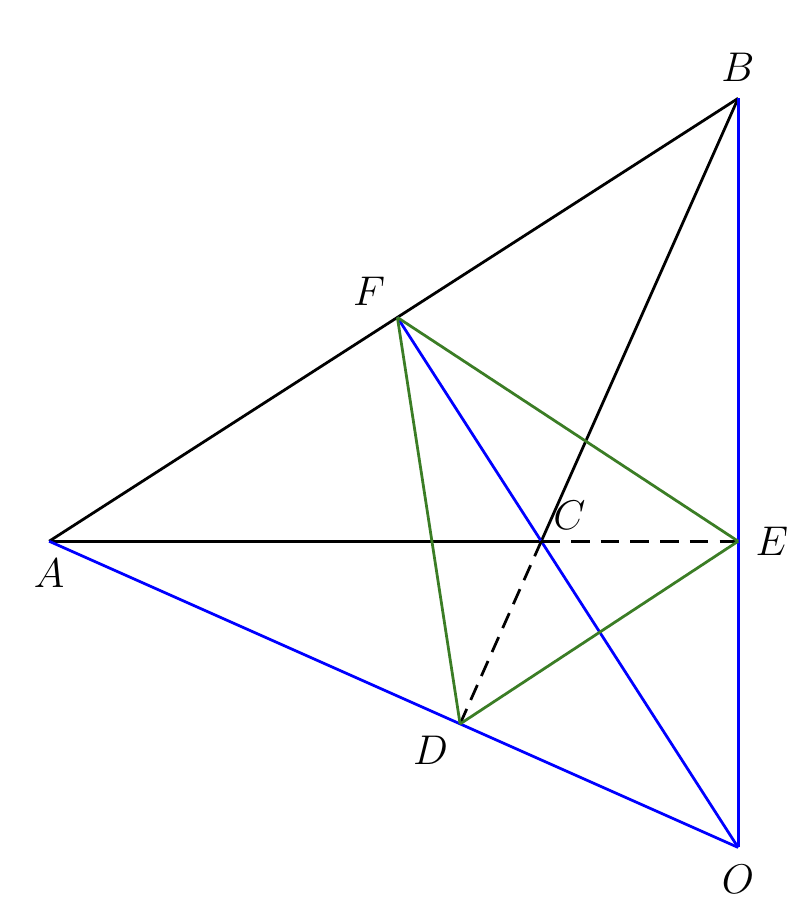}}\\
b)
\end{tabular}
\caption{Alturas en un tri'angulo oblicu'angulo cuando: a) el tri'angulo es acut'angulo, b) el tri'angulo es obtus'angulo.}\label{F17}
\end{figure}

La altura de un tri'angulo oblicu'angulo es el segmento formado por una recta que pasa por un v'ertice y es perpendicular a la recta que contiene a su lado opuesto. Consecuencia de esta definici'on es que un tri'angulo cualquiera tendr'a tres alturas. Una propiedad muy conocida de estas es que se intersectan en un punto denominado ortocentro (ver la figura \ref{F17}).

Sea $h_a=\seg{AD}$, $h_b=\seg{BE}$ y $h_c=\seg{CF}$ las longitudes de las alturas del tri'angulo $ABC$. De los tri'angulos rect'angulos $ACD$, $ABD$, $ABE$, $BCE$, $ACF$, $BCF$, se tiene que 
\begin{align}
h_a&=b\sen{C}=c\sen{B}\label{E201}\\
h_b&=c\sen{A}=a\sen{C}\label{E202}\\
h_c&=a\sen{B}=b\sen{A}\label{E203}
\end{align}
Utilizando la identidad $\sen{\al}=2\sen{\frac{\al}{2}}\cos{\frac{\al}{2}}$ en los miembros medios de las ecuaciones anteriores y utilizando las f'ormulas del 'angulo medio, se obtienen expresiones de las alturas en funci'on de los lados del tri'angulo $ABC$. Estas son:
\begin{align}
h_a=b\sen{C}&=2b\sqrt{\frac{(p-a)(p-b)}{ab}}\sqrt{\frac{p(p-c)}{ab}}\notag\\
&=\frac{2}{a}\sqrt{p(p-a)(p-b)(p-c)}\label{E204}\\
h_b=c\sen{A}&=2c\sqrt{\frac{(p-b)(p-c)}{bc}}\sqrt{\frac{p(p-a)}{bc}}\notag\\
&=\frac{2}{b}\sqrt{p(p-a)(p-b)(p-c)}\label{E205}\\
h_c=a\sen{B}&=2a\sqrt{\frac{(p-a)(p-c)}{ac}}\sqrt{\frac{p(p-b)}{ac}}\notag\\
&=\frac{2}{c}\sqrt{p(p-a)(p-b)(p-c)}\label{E206}
\end{align}
Si se suman los inversos de las alturas:
\begin{equation}\label{E207}
\frac{1}{h_a}+\frac{1}{h_b}+\frac{1}{h_c}=\frac{a/2+b/2+c/2}{\sqrt{p(p-a)(p-b)(p-c)}}=\frac{1}{r}
\end{equation}
Combinando esta ecuaci'on con la ecuaci'on \eqref{E139} se tiene que
\begin{equation}\label{E208}
\frac{1}{h_a}+\frac{1}{h_b}+\frac{1}{h_c}=\frac{1}{r_a}+\frac{1}{r_b}+\frac{1}{r_c}
\end{equation}
Esta relaci'on, en palabras, dice que el promedio arm'onico de las alturas de un tri'angulo oblicu'angulo es igual al promedio arm'onico de los radios de sus circunferencias excritas; promedios que son iguales al radio de su circunferencia inscrita.

Si se unen los pies de las alturas $D$, $E$ y $F$ del tri'angulo, se obtiene un tri'angulo conocido como pedal, el cual puede verse en las figuras \ref{F17}-a y \ref{F17}-b en verde para los tri'angulos acut'angulo y obtus'angulo, respectivamente. 

Para resolver este tri'angulo, es necesario caracterizar los tri'angulos $AEF$, $BDF$ y $CDE$. Se comienza con el tri'angulo $CDE$ de la figura \ref{F17}-a, del cual se conoce los lados $\seg{CD}=b\cos{C}$ y $\seg{CE}=a\cos{C}$ y el 'angulo $\angle DCE=C$ (con respecto al tri'angulo correspondiente en la figura \ref{F17}-b, $\seg{CD}=b\cos{(180^\circ-C)}=-b\cos{C}$ y $\seg{CE}=a\cos{(180^\circ-C)}=-a\cos{C}$; segmentos que son positivos, pues $C$ es obtuso y su coseno ser'a negativo). Aplicando el teorema del coseno, se tiene que
\begin{align*}
\seg{DE}^2&=\seg{CD}^2+\seg{CE}^2-2\seg{CD}\cdot \seg{CE}\cos{C}=b^2\cos^2{C}+a^2\cos^2{C}-2ab\cos^3{C}\\
&=(b^2+a^2-2ab\cos{C})\cos^2{C}
\end{align*}
Del tri'angulo $ABC$, se sabe que $c^2=b^2+a^2-2ab\cos{C}$, luego
\begin{equation}\label{E209}
\seg{DE}=c|\cos{C}|
\end{equation}
Los 'angulos internos restantes podr'ian determinarse mediante el teorema del seno, sin embargo, si se observa que los tri'angulos $ABC$ y $CDE$ tienen un 'angulo igual ($\angle BCA=\angle DCE$) y que el cociente entre los lados $\seg{CD}$ y $\seg{AC}$ es igual al cociente de los lados $\seg{CE}$ y $\seg{BC}$ (el cual es $\cos{C}$), entonces ellos son semejantes, y por consiguiente, $\angle DEC=B$ y $\angle EDC=A$. Esta observaci'on se aplica igualmente a los tri'angulos correspondientes de la figura \ref{F17}-b, solo que la igualdad de los 'angulos se deduce porque son opuestos al v'ertice y el cociente entre los lados es $\cos{(180^\circ-C)}$. En resumen: Para el tri'angulo $CDE$:
\begin{align}
\seg{DE}&=c|\cos{C}|\label{E210}\\
\angle DEC&=B\label{E211}\\
\angle EDC&=A\label{E212}
\end{align}
Si se sigue un an'alisis an'alogo con los tri'angulos restantes, se tendr'a lo siguiente: Para el tri'angulo $BDF$:
\begin{align}
\seg{DF}&=b|\cos{B}|\label{E213}\\
\angle DFB&=C\label{E214}\\
\angle FDB&=A\label{E215}
\end{align}
y para el tri'angulo $AEF$:
\begin{align}
\seg{EF}&=a|\cos{A}|\label{E216}\\
\angle FEA&=B\label{E217}\\
\angle EFA&=C\label{E218}
\end{align}
Esto aplica para los tri'angulos correspondientes de la construcci'on del tri'anguo obtus'angulo. Del tri'angulo $DEF$ se conoce sus lados, dados por las ecuaciones \eqref{E213}, \eqref{E216} y \eqref{E219}. Sus 'angulos internos, sin embargo, deben calcularse con cuidado. 

Comenzando con el 'angulo $\angle EFD$, en la construcci'on del tri'angulo acut'angulo:
\[
\angle EFD=180^\circ-(\angle EFA+\angle DFB)=180^\circ-2C
\]
y en la construcci'on del tri'angulo obtus'angulo:
\begin{align*}
\angle EFD&=180^\circ-(\angle DFA+\angle EFB)=180^\circ-(180^\circ-\angle EFA+180^\circ-\angle DFB)
\\&=2C-180^\circ
\end{align*}
este 'angulo puede expresarse en una sola ecuaci'on:
\begin{equation}\label{E219}
\angle EFD=|180^\circ-2C|
\end{equation}
Los 'angulos $\angle FDE$ y $DEF$, mediante un an'alisis similar son iguales a 
\begin{align}
\angle FDE&=|180^\circ-2A|\label{E220}\\
\angle DEF&=|180^\circ-2B|\label{E221}
\end{align}
Esta secci'on finaliza con una discusi'on sobre c'omo resolver un tri'angulo obligu'angulo dadas sus alturas. Para ello, se expresa el 'area $S$ del tri'angulo en funci'on de sus alturas, es decir:
\begin{equation}\label{E222}
2S=ah_a=bh_b=ch_c
\end{equation}
Si se reescribe este conjunto de igualdades de la siguiente manera:
\begin{equation}\label{E223}
2S=\cfrac{a}{\cfrac{1}{h_a}}=\cfrac{b}{\cfrac{1}{h_b}}=\cfrac{c}{\cfrac{1}{h_c}}
\end{equation}
entonces, de las propiedades de las proporciones, se tiene que
\begin{equation}\label{E224}
2S=\cfrac{2p}{\cfrac{1}{h_a}+\cfrac{1}{h_b}+\cfrac{1}{h_c}}
=\cfrac{p}{\cfrac{1}{2}\GP{\cfrac{1}{h_a}+\cfrac{1}{h_b}+\cfrac{1}{h_c}}}
\end{equation}
Con las ecuaciones \eqref{E223} y \eqref{E224}, se puede construir las siguientes proporciones:
\begin{align}
2S&=\cfrac{p-a}{\cfrac{1}{2}\GP{\cfrac{1}{h_b}+\cfrac{1}{h_c}-\cfrac{1}{h_a}}}\label{E225}\\
2S&=\cfrac{p-b}{\cfrac{1}{2}\GP{\cfrac{1}{h_a}+\cfrac{1}{h_c}-\cfrac{1}{h_b}}}\label{E226}\\
2S&=\cfrac{p-c}{\cfrac{1}{2}\GP{\cfrac{1}{h_a}+\cfrac{1}{h_b}-\cfrac{1}{h_c}}}\label{E227}
\end{align}
Con ellas, se puede expresar tanto el semiper'imetro como su diferencia con los lados en funci'on del 'area y de las alturas, es decir:
\begin{align}
p&=S\GP{\cfrac{1}{h_a}+\cfrac{1}{h_b}+\cfrac{1}{h_c}}\label{E228}\\
p-a&=S\GP{\cfrac{1}{h_b}+\cfrac{1}{h_c}-\cfrac{1}{h_a}}\label{E229}\\
p-b&=S\GP{\cfrac{1}{h_a}+\cfrac{1}{h_c}-\cfrac{1}{h_b}}\label{E230}\\
p-c&=S\GP{\cfrac{1}{h_a}+\cfrac{1}{h_b}-\cfrac{1}{h_c}}\label{E231}
\end{align}
Ahora, empleando las f'ormulas de la tangente del 'angulo medio, se puede calcular los 'angulos internos:
\begin{align}
\tg{\frac{A}{2}}&=\sqrt{\frac{(p-b)(p-c)}{p(p-a)}}=\sqrt{\cfrac{\GP{\cfrac{1}{h_a}+\cfrac{1}{h_c}-\cfrac{1}{h_b}}\GP{\cfrac{1}{h_a}+\cfrac{1}{h_b}-\cfrac{1}{h_c}}}{\GP{\cfrac{1}{h_a}+\cfrac{1}{h_b}+\cfrac{1}{h_c}}\GP{\cfrac{1}{h_b}+\cfrac{1}{h_c}-\cfrac{1}{h_a}}}}\label{E232}\\
\tg{\frac{B}{2}}&=\sqrt{\frac{(p-a)(p-c)}{p(p-b)}}=\sqrt{\cfrac{\GP{\cfrac{1}{h_b}+\cfrac{1}{h_c}-\cfrac{1}{h_a}}\GP{\cfrac{1}{h_a}+\cfrac{1}{h_b}-\cfrac{1}{h_c}}}{\GP{\cfrac{1}{h_a}+\cfrac{1}{h_b}+\cfrac{1}{h_c}}\GP{\cfrac{1}{h_a}+\cfrac{1}{h_c}-\cfrac{1}{h_b}}}}\label{E233}\\
\tg{\frac{C}{2}}&=\sqrt{\frac{(p-a)(p-b)}{p(p-c)}}=\sqrt{\cfrac{\GP{\cfrac{1}{h_b}+\cfrac{1}{h_c}-\cfrac{1}{h_a}}\GP{\cfrac{1}{h_a}+\cfrac{1}{h_c}-\cfrac{1}{h_b}}}{\GP{\cfrac{1}{h_a}+\cfrac{1}{h_b}+\cfrac{1}{h_c}}\GP{\cfrac{1}{h_a}+\cfrac{1}{h_b}-\cfrac{1}{h_c}}}}\label{E234}
\end{align}
Estas ecuaciones pueden simplificarse si se utiliza la ecuaci'on \eqref{E207}. Al hacerlo, estas se convierten en:
\begin{align}
\tg{\frac{A}{2}}&=\sqrt{\cfrac{\GP{\cfrac{1}{r}-\cfrac{2}{h_b}}\GP{\cfrac{1}{r}-\cfrac{2}{h_c}}}{\cfrac{1}{r}\GP{\cfrac{1}{r}-\cfrac{2}{h_a}}}}=
\sqrt{\frac{h_a(h_b-2r)(h_c-2r)}{h_bh_c(h_a-2r)}}\label{E235}\\
\tg{\frac{B}{2}}&=\sqrt{\cfrac{\GP{\cfrac{1}{r}-\cfrac{2}{h_a}}\GP{\cfrac{1}{r}-\cfrac{2}{h_c}}}{\cfrac{1}{r}\GP{\cfrac{1}{r}-\cfrac{2}{h_b}}}}=
\sqrt{\frac{h_b(h_a-2r)(h_c-2r)}{h_ah_c(h_b-2r)}}\label{E236}\\
\tg{\frac{C}{2}}&=\sqrt{\cfrac{\GP{\cfrac{1}{r}-\cfrac{2}{h_a}}\GP{\cfrac{1}{r}-\cfrac{2}{h_b}}}{\cfrac{1}{r}\GP{\cfrac{1}{r}-\cfrac{2}{h_c}}}}=
\sqrt{\frac{h_c(h_a-2r)(h_b-2r)}{h_ah_b(h_c-2r)}}\label{E237}
\end{align}
Las ecuaciones anteriores, adem'as de expresar los 'angulos internos en funci'on de las alturas, da una condici'on para saber cu'ando puede resolverse un tri'angulo si se conocen sus alturas, las cuales deben ser mayores que el di'ametro de la circunferencia inscrita a dicho tri'angulo, o si se prefiere expresarla de otra manera, las alturas deben ser mayores al doble de su media arm'onica. Sobra decir que, una vez calculados los 'angulos, debe verificarse que su suma sea de $180^\circ$. Si esto se cumple, los lados se calculan mediante las ecuaciones \eqref{E201}--\eqref{E203}:
\begin{align}
a&=\frac{h_b}{\sen{C}}=\frac{h_c}{\sen{B}}\label{E238}\\
b&=\frac{h_c}{\sen{A}}=\frac{h_a}{\sen{C}}\label{E239}\\
c&=\frac{h_a}{\sen{B}}=\frac{h_b}{\sen{A}}\label{E240}
\end{align}
La soluci'on obtenida se verifica primero, examinando que el semiper'imetro sea mayor a los lados; y segundo, determinar si se satisface una de las relaciones de la proyecci'on o de Mollweide.
%\begin{equation}
%\end{equation}
\addcontentsline{toc}{section}{Bibliograf'ia}
\renewcommand{\refname}{Bibliograf\'ia}

\end{document}